\documentclass[titlepage,final,11pt]{article}
\usepackage{theorem}
\usepackage{enumerate}
\usepackage{array}
\usepackage[reqno]{amsmath}
\usepackage{amssymb}
\usepackage{mathtools}
\usepackage{latexsym}
\usepackage{makeidx}
\usepackage{fancybox}
\input epsf.sty
\usepackage[dvips]{graphicx,color}
\usepackage{showlabels}
\usepackage{subcaption}
\usepackage{float}

\newcommand{\drawing}[4]{
	\begin{figure}[!hbt]
	\begin{center}
	\leavevmode
	\epsfxsize=#2
	\epsfbox{#1}
	\caption{\small #3}
	\label{#4}
	\end{center}
	\end{figure}}

\theoremstyle{change}
\newtheorem{proclaim}{PROCLAIM}[section]
\newtheorem{theorem}[proclaim]{Theorem}

\newtheorem{proposition}[proclaim]{Proposition}
\newtheorem{corollary}[proclaim]{Corollary}

\numberwithin{equation}{section}

\outer\def\proclaim #1. #2\par{\medbreak \noindent{\bf#1.\enspace}{\sl#2}\par
  \ifdim\lastskip<\medskipamount
  \removelastskip\penalty55\medskip\fi}
\def\state #1. { \noindent{\bf#1.\enspace}}
\def\algo #1. { \noindent{\bf#1.\enspace}}

\DeclareMathOperator{\nt}{int}

\newcommand{\comp}{\,{\raise 1pt \hbox{$\scriptstyle\circ$}}\,}

\newcommand{\reals}{\mathbb{R}}

\newcommand{\natnums}{{{\rm l} \kern -.13em {\rm N} }}
\newcommand{\nats}{\mathbb{N}}
\newcommand{\snats}{{I\kern -.29em N}}
\newcommand{\rats}{{Q\kern -.64em \raise 1pt \hbox{$\scriptstyle |$}\;\,}}
\newcommand{\srats}
	{{Q\kern -.56em \raise 1.2pt \hbox{$\scriptscriptstyle /$}\,}}
\newcommand{\ints}{Z\kern -.46em Z}
\newcommand{\ball}{\mathbb{B}}
\newcommand{\pluss}{\hskip1pt \raise1pt\vbox{\hrule width6pt \vskip1pt \hrule
                    width6pt} \kern-4pt{\lower1pt\hbox{\vrule height6pt
		    \kern1pt\vrule height6pt}}\hskip5pt}
\newcommand{\eop}
	{\hfill{$\vcenter{\hrule height1pt \hbox{\vrule width1pt height5pt
   	 \kern5pt \vrule width1pt} \hrule height1pt}$} \medskip}

\newcommand{\setd}{{ d \kern -.15em l}}
\newcommand{\hatsetd}{ d \hat{\kern -.15em l }}

\renewcommand{\epsilon}{\varepsilon}
\renewcommand{\phi}{\varphi}

\newcommand{\tto}{\;{\lower 1pt \hbox{$\rightarrow$}}\kern -12pt
           \hbox{\raise 2.5pt \hbox{$\rightarrow$}}\;}
\newcommand{\overto}[1]{\,{\raise 0pt\hbox{$\rightarrow$}}\kern -9pt
     \hbox{\lower 3pt \hbox{$\scriptscriptstyle#1$}}\hskip6pt}
\newcommand{\underto}[1]{\,{\lower 1pt\hbox{$\rightarrow$}}\kern -9pt
     \hbox{\raise 4pt \hbox{$\,\scriptscriptstyle#1$}}\hskip7pt}
\newcommand{\bigoverto}[1]{{\raise 0pt\hbox{$\,\longrightarrow$}}\kern -16pt
     \hbox{\lower 3pt \hbox{$\scriptscriptstyle#1$}}\hskip4pt}
\newcommand{\bigunderto}[1]{\,{\lower 1pt\hbox{$\longrightarrow$}}\kern -16pt
     \hbox{\raise 4pt \hbox{$\,\scriptscriptstyle#1$}}\hskip6pt}
\newcommand{\bigbigto}[2]{\,{\raise 0pt\hbox{$\,\longrightarrow$}}\kern -16pt
     \hbox{\lower 3pt \hbox{$\scriptscriptstyle#2$}}\kern -10pt
     \hbox{\raise 4pt \hbox{$\,\scriptscriptstyle#1$}}\hskip7pt}
\newcommand{\downto}{{\raise 1pt \hbox{$\scriptscriptstyle \,\searrow\,$}}}
\newcommand{\upto}{{\raise 1pt \hbox{$\scriptscriptstyle \,\nearrow\,$}}}

\newcommand{\notimply}
	{\quad\hbox{$\Longrightarrow \kern -14pt {/}$}\hskip6pt\quad}

\newcommand{\lto}{\,{\lower 1pt\hbox{$\rightarrow$}}\kern -10pt
     \hbox{\raise 4pt \hbox{$\, \scriptstyle l$}}\hskip7pt}
\newcommand{\eto}{\,{\lower 1pt\hbox{$\rightarrow$}}\kern -11pt
     \hbox{\raise 4pt \hbox{$\, \scriptstyle e$}}\hskip7pt}
\newcommand{\hto}{\,{\lower 1pt\hbox{$\rightarrow$}}\kern -11pt
     \hbox{\raise 4pt \hbox{$\, \scriptstyle h$}}\hskip7pt}
\newcommand{\pto}{\,{\lower 1pt\hbox{$\rightarrow$}}\kern -11pt
     \hbox{\raise 4.5pt \hbox{$\, \scriptstyle p$}}\hskip7pt}
\newcommand{\cto}{\,{\lower 1pt\hbox{$\rightarrow$}}\kern -11pt
     \hbox{\raise 4pt \hbox{$\, \scriptstyle c$}}\hskip7pt}
\newcommand{\gto}{\,{\lower 1pt\hbox{$\rightarrow$}}\kern -11pt
     \hbox{\raise 4.5pt \hbox{$\, \scriptstyle g$}}\hskip7pt}
\newcommand{\sto}{\,{\lower 1pt\hbox{$\rightarrow$}}\kern -11pt
     \hbox{\raise 4pt \hbox{$\, \scriptstyle s$}}\hskip7pt}
\newcommand{\awto}{\,{\lower 1pt\hbox{$\rightarrow$}}\kern -15pt
     \hbox{\raise 4pt \hbox{$\, \scriptstyle aw$}}\hskip7pt}
\def\Nto{\,{\raise 1pt\hbox{$\rightarrow$}}\kern -13pt
     \hbox{\lower 3pt \hbox{$\, \scriptstyle N$}}\hskip7pt}
\def\Cto{\,{\raise 1pt\hbox{$\rightarrow$}}\kern -14pt
     \hbox{\lower 3pt \hbox{$\, \scriptstyle C$}}\hskip7pt}
\def\fto{\,{\raise 1pt\hbox{$\rightarrow$}}\kern -14pt
     \hbox{\lower 3pt \hbox{$\, \scriptstyle f$}}\hskip7pt}

\newcommand{\low}[1]{{\lower1pt \hbox{$\scriptstyle #1$}}}
\newcommand{\loww}[1]{{\lower2pt \hbox{$\scriptstyle #1$}}}
\newcommand{\high}[1]{{\raise1pt \hbox{$\scriptstyle #1$}}}

\newcommand{\cF}{{\cal F}}
\newcommand{\cG}{{\cal G}}

\newcommand{\cL}{{\cal L}}

\newcommand{\cR}{{\cal R}}

\newcommand{\cZ}{{\cal Z}}

\newcommand{\nliminf}{\mathop{\rm liminf}\nolimits}

\newcommand{\nlimsup}{\mathop{\rm limsup}\nolimits}

\newcommand{\nnmin}{\mathop{\rm minimize}}

\newcommand{\nargmin}{\mathop{\rm argmin}\nolimits}

\newcommand{\var}{{\mathop{\rm var}\nolimits}}
\newcommand{\std}{\hspace{0.04cm}{\mathop{\rm std}\nolimits}\hspace{-0.04cm}}
\newcommand{\prob}{{\mathop{\rm prob}\nolimits}}

\newcommand{\bfxi}{\mbox{\boldmath $\xi$}}
\newcommand{\bfeta}{\mbox{\boldmath $\eta$}}

\newcommand{\bfL}{\mbox{\boldmath $L$}}

\newcommand{\bfy}{\mbox{\boldmath $y$}}

\newcommand{\bfx}{\mbox{\boldmath $x$}}

\newcommand{\lwdy}[2]{\mathrel{\mathop
        {\raisebox{0.1ex}{\null$#1$}}{\hbox{\kern -1.0em
	{\raisebox{-0.8ex}{$\scriptstyle{\;\to #2}$}}}}}}
\newcommand{\lwwdy}[2]{\mathrel{\mathop
        {\raisebox{0.2ex}{\null$#1$}}{\hbox{\kern -1.0em
	{\raisebox{-1.1ex}{$\scriptstyle{\;\to #2}$}}}}}}
\newcommand{\slwdy}[2]{\scriptsize{{\mathrel{\mathop
        {\raisebox{0.1ex}{\null$#1$}}{\hbox{\kern -1.0em
	{\raisebox{-0.8ex}{$\scriptstyle{\;\to #2}$}}}}}}}}
\newcommand{\slwwdy}[2]{\scriptsize{{\mathrel{\mathop
        {\raisebox{0.2ex}{\null$#1$}}{\hbox{\kern -1.0em
	{\raisebox{-1.1ex}{$\scriptstyle{\;\to #2}$}}}}}}}}

\definecolor{lightgray}{gray}{0.75}
\definecolor{myred}{rgb}{0.55,0,0}
\definecolor{myblue}{rgb}{0,0,0.5} 
\definecolor{mygreen}{rgb}{0,0.5,0} 
\definecolor{purple}{rgb}{0.5,0,0.5} 
\definecolor{turq}{rgb}{0,0.805,0.816} 
\definecolor{maroon}{rgb}{0.51,0,0}
\definecolor{MAROON}{rgb}{0.51,0,0}
\definecolor{redor}{rgb}{0.78,0.078,0.078}
\definecolor{dgreen}{rgb}{0,0.3,0}

\newcommand{\Ex}{\mathbb{E}}

\newcommand{\bcdot}{\,{\raise .2ex \hbox{$\centerdot$}}\,}

\newcommand{\diam}{\textup{\textrm{diam}}}

\newcommand{\bbH}{\mathbb{H}}

\usepackage{multirow}
\usepackage{graphicx}
\usepackage{epsfig}
\usepackage{adjustbox}
\usepackage{rotating}
\usepackage{float}
\usepackage{color}
\usepackage{makecell}
\usepackage{xcolor}
\usepackage{cases}
\usepackage{empheq}
\usepackage{placeins}
\usepackage{qtree}
\usepackage{alphalph}
\usepackage{tabulary}
\usepackage{array}
\usepackage{paralist}
\usepackage{graphics}
\usepackage{caption}
\usepackage{subcaption}
\usepackage{cases}
\usepackage{nccmath}
\usepackage{tabularx,booktabs}
\usepackage{comment}

\begin{document}


\begin{center}
\begin{large}
{\bf Risk-Adaptive Local Decision Rules}
\smallskip
\end{large}
\vglue 0.7truecm
\begin{tabular}{lcl}
  \begin{large}
  {\sl Johannes O. Royset }
  \end{large} & \ \ {\phantom{\&}} \ \ &
   \begin{large} {\sl Miguel A. Lejeune
   				  } \end{large} \\
  \\
  Daniel J. Epstein Department of &&   Department of Decision Sciences \\
  Industrial \& Systems Engineering && School of Business \\
  University of Southern California  && George Washington University \\
  royset@usc.edu && mlejeune@gwu.edu
\end{tabular}

\vskip 0.2truecm

\end{center}

\vskip 0.5truecm

\noindent {\bf Abstract}. For parameterized mixed-binary optimization problems, we construct local decision rules that prescribe near-optimal courses of action across a set of parameter values. The decision rules stem from solving risk-adaptive training problems over classes of continuous, possibly nonlinear mappings. In asymptotic and nonasymptotic analysis, we establish that the decision rules prescribe near-optimal decisions locally for the actual problems, without relying on linearity, convexity, or smoothness. The development also accounts for practically important aspects such as inexact function evaluations, solution tolerances in training problems, regularization, and reformulations to solver-friendly models. The decision rules also furnish a means to carry out sensitivity and stability analysis for broad classes of parameterized optimization problems. We develop a decomposition algorithm for solving the resulting training problems and demonstrate its ability to generate quality decision rules on a nonlinear binary optimization model from search theory.

\vskip 0.2truecm

\halign{&\vtop{\parindent=0pt
   \hangindent2.5em\strut#\strut}\cr
                         \cr

{\bf Date}:\quad \ \today \cr}

\baselineskip=15pt

\section{Introduction}\label{sec:intro}

A decision rule prescribes courses of action as alternatives to those obtained by solving an optimization problem. It is usually easy to interpret and explain, and thus suitable for decisions requiring transparency and buy-in from stakeholders and the general public, which are increasingly important as discussed in \cite{GoerigkHartisch.23}. Decision rules allow for quicker decision making than relying on the solution of large-scale optimization problems. A decision rule also tends to prescribe similar decisions for similar situations, while poorly formulated optimization problems might have solutions that are overly sensitivity to small changes to parameters \cite{Royset.20b,Royset.21}. We refer to the review articles \cite{DelageIancu.15,GeorghiouKuhnWiesemann.19,YanikogluGorissenHertog.19} for a general introduction. 

In this paper, we develop a framework for constructing local decision rules for parameterized problems with mixed-binary decision variables. We use risk measures to adapt the resulting decision rules to the level of conservativeness deemed necessary by the analyst. The rules turn out to be locally consistent relative to an optimal decision for a given parameter vector, without imposing structural assumptions such as linearity, convexity, and/or smoothness about the underlying problem; the decision rules only need to be continuous mappings. Thus, for broad classes of problems, the framework produces a local decision rule that is ``nearly'' as good as an optimal solution tailored to a specific parameter vector. In a nonasymptotic analysis, we quantify the level of suboptimality for the decision rules as they prescribe courses of actions for different parameter vectors, again without relying on linearity, convexity, and/or smoothness. The development accounts for practically important aspects such as inexact computation of functions, solution tolerances in training problems, regularization, reformulations to solver-friendly models, and decomposition for scalable computations. The framework also provides a means to carry out sensitivity and stability analysis of parameterized optimization problems. 

For a parameter vector $\xi\in \reals^r$, we consider the parameterized mixed-binary optimization problem
\begin{equation*}\label{eqn:actualproblem}
\mbox{(AP)}(\xi)~~~~~\nnmin_{(x, y) \in C} \, \psi_0(\xi,x,y) ~~\mbox{ subject to } ~~\psi_k(\xi,x,y) \leq 0, ~~k=1, \dots, q,
\end{equation*}
referred to as the {\em actual problem}, where $C\subset \reals^n \times \{0,1\}^m$ is a nonempty and closed set and the functions $\psi_k:\reals^r\times \reals^n\times \reals^m\to \reals$, $k=0,1, \dots, q$, are continuous throughout the paper. Thus, $x$ is a vector of real variables and $y$ is a vector of binary variables. While one might be able to solve the problem for each $\xi$, we seek to develop a decision rule that prescribes a reasonably good mixed-binary vector $(x,y)\in C$ for each $\xi$ of interest. We view a decision rule as a pair of mappings $(F,G)$, with $F:\reals^r\to \reals^n$ and $G:\reals^r\to \reals^m$. For any $\xi\in \reals^r$, $F(\xi)\in \reals^n$ prescribes a decision $x$ in the actual problem (AP)$(\xi)$. The binary vector $y$ is similarly defined but leverages the Heaviside function. For any positive integer $s$, let $\bbH:\reals^s\to \reals^s$ be the multi-variate Heaviside function, i.e., for $v = (v_1, \dots, v_s)\in \reals^s$,
\[
\bbH(v) = (w_1, \dots, w_s), \mbox{ with } w_i = \begin{cases}
0 & \mbox{ if } v_i \leq 0\\
1 & \mbox{ otherwise.}
\end{cases}
\]
For each $\xi$, the mapping $G$ defines an $m$-dimensional binary vector $\bbH(G(\xi))$, which in turn prescribes a decision $y$ in (AP)$(\xi)$; see \cite{BertsimasGeorghiou.18} for similar usage of the Heaviside function to define decision rules.

Given a class $\cF$ of mappings from $\reals^r$ to $\reals^n$ and a class $\cG$ of mappings from $\reals^r$ to $\reals^m$, we seek to determine $F\in \cF$ and $G\in \cG$ such that $(F(\xi), \bbH(G(\xi)))$ is a good solution (in some sense) for (AP)$(\xi)$ and that this holds for a collection of parameter vectors $\xi$ of interest. The process of finding $F$ and $G$ is a {\em training problem}. The formulation and solution of such training problems are main challenges discussed in this paper. We consider broad classes of mappings and essentially only require that they are continuous. The mappings can be given by neural nets or by polynomial, piecewise affine, or affine functions. 

The appeal of decision rules is long-standing and their systematic study extends at least back to \cite{CharnesCooperSymonds.58}; see also \cite{GarstkaWets.74}. The seminal paper \cite{BentalGoryashkoGuslitzerNemirovski.04} reinvigorated the field and spurred further developments in several directions over the last two decades under the name adjustable robust optimization. While classical robust optimization (see, e.g., \cite{BentalElghaouiNemirovski.09}) determines a decision $(x,y)$ that is satisfactory in AP($\xi$) for a set of parameter vectors $\xi$, adjustable robust optimization takes a cue from two-stage stochastic optimization (see, e.g., \cite{Wets.74} and \cite[Chapter 3]{primer}) and adapts some of the decision variables to specific values of $\xi$. One seeks a decision rule for those variables, but this is generally computationally difficult to achieve \cite[Section 14.2]{BentalElghaouiNemirovski.09}. Nevertheless, there have been successes starting with affine decision rules for linear programs subject to parameter uncertainty \cite{BentalGoryashkoGuslitzerNemirovski.04,ChenSimSunZhang.08}. Generally, one can expect affine decision rules to prescribe suboptimal decisions as the minimizers of AP($\xi$) tend to vary nonlinearly in $\xi$ but proof of optimality exists in specific cases; see, e.g., classical results from control theory of linear systems as well as the more recent papers \cite{BertsimasIancuParrilo.10,GounarisWiesemannFloudas.13,ArdestaniJaafariDelage.16,SimchiLeviTrichakisZhang.19,YanikogluGorissenHertog.19,GeorghiouTsoukalasWiesemann.21}. Various guarantees about the optimality gap for affine decision rules appear in \cite{KuhnWiesemannGeorghiou.11,BertsimasGoyal.13,BertsimasBidkhori.15,ElhousniGoyal.17,ElHousniGoyal.21}. A constant decision rule (as in classical robust optimization) might also suffice; see \cite{BertsimasGoyalLu.15,AwasthiGoyalLu.19} for a quantification of the optimality gap in such cases. 

As introduced in \cite{BertsimasCaramanis.10} (see also \cite{BertsimasCaramanis.07} and its reference to a presently unavailable 2007 technical report), the approach of finite adaptivity considers piecewise constant decision rules; see also \cite{VayanosKuhnRustem.11,HanasusantoKuhnWiesemann.15,BertsimasGeorghiou.18} and \cite{PostekHertog.16,BertsimasDunning.16} for means to iteratively partition the set of $\xi$-values of interest. The presence of binary variables (such as $y$ in AP($\xi$)) or more generally integer variables represents a major challenge; see \cite{BertsimasCaramanis.07} for an early effort based on integer matrices and vectors in an affine decision rule combined with ceiling fractional values of $\xi$. Finite adaptivity, while already used to tackle integer variables in \cite{BertsimasCaramanis.10}, is further refined for that purpose in \cite{HanasusantoKuhnWiesemann.15}, which identifies the number of candidate decision vectors needed in specific cases and proposes mixed-integer linear programming formulations; see also \cite{SubramanyamGounarisWiesemann.20}. Other reformulations that rely on duality also lead to mixed-integer linear programs as described in \cite{GeorghiouWiesemannKuhn.15}. 

While various lifting techniques extend the range of affine decision rules (see, e.g., \cite{ChenZhang.09,GeorghiouWiesemannKuhn.15}), it is apparent that more flexible decision rules are needed in many applications. Piecewise linear decision rules appear in \cite{ChenSimSunZhang.08,ChenZhang.09,GohSim.10,BertsimasGeorghiou.15,GeorghiouWiesemannKuhn.15} and they may also leverage neural networks \cite{RahalLiPapageorgiou.22}. Polynomial decision rules are studied in \cite{BampouKuhn.11,BertsimasIancuParrilo.11,GeorghiouWiesemannKuhn.15}; see Section 7 of \cite{YanikogluGorissenHertog.19} for a review. 

With few exceptions (cf. \cite{TakedaTaguchiTutuncu.08}), the existing literature focuses on linear optimization problems (possibly with integer variables) and this is well-motivated by computational considerations. The contributions of the present paper stem from its focus on arbitrary real-valued continuous functions $\psi_0, \psi_1, \dots, \psi_q$  and an arbitrary closed set $C$ in (AP)($\xi$). In such a general setting, we show that novel training problems defined using measures of risk and leveraging a finite training set $\{\xi^1, \dots, \xi^s\in \reals^r\}$ produce decision rules that are locally consistent in general and, in the absence of the constraints $\psi_k(\xi,x,y) \leq 0$, $k=1, \dots, q$, also possess nonasymptotic suboptimality guarantees. The training problems can be solved with a nonzero tolerance and can involve inaccurate evaluations of the functions $\psi_0, \psi_1, \dots, \psi_q$. In the absence of  $\psi_k(\xi,x,y) \leq 0$, $k=1, \dots, q$, we show that the obtained decision rules are $\epsilon$-suboptimal relative to (AP)($\xi$) across a compact subset of parameter vectors $\xi$, with $\epsilon$ being the sum of two times the inaccuracy in function evaluation, the optimality tolerance in the training problem, and two times the product of the diameter of the compact subset and a Lipschitz modulus of $\psi_0$; see Theorem \ref{thm:error}. The decision rules are nearly arbitrary continuous mappings and thus we cover affine ones as well as many other possibilities within the same analysis. The choice of measures of risk in the construction of training problems is also nearly arbitrary and certainly permits expected value, worst-case, quantiles (a.k.a. VaR), and superquantiles (a.k.a. CVaR, AVaR, expected shortfall). Thus, we can adjust the level of conservativeness in a training problem to match the need in a particular setting. The choice of expected value amounts to asking for a decision rule that is good {\em on average} across the training set. This choice does not rule out the possibility of poor performance by the decision rule at {\em some} training points. A worst-case risk measure results in a decision rule that is reasonably good across the whole training set. The use of a superquantile risk measure allows us to slide between these two extremes; recall that the $\alpha$-superquantile of a random variable is the average of the worst $(1-\alpha)100\%$ outcomes. In effect, users can draw from a wide array of risk measures to capture their preferences. Empirical and theoretical evidence from supervised learning and elsewhere indicate that a risk-averse criterion might result in better generalization than when relying on expected values \cite{LevyCarmonDuchiSidford.20,LaguelPillutlaMalickHarchaoui.21b,LaguelPillutlaMalickHarchaoui.21,GotohKimLim.21}. These indications motivate our consideration of general risk measures and, indeed, our empirical study confirms that superquantile risk measures typically provide better decision rules than expected values when tested on $\xi$-values outside the training set. Still, we shy away from recommending a {\em particular} risk measure as the choice necessarily will depend on the specific setting. As with tuning of ``hyper-parameters'' in machine learning, a verification data set can be brought in to help with the choice.

The paper also develops a decomposition algorithm for solving typical instances of the training problems and derives upper and lower bounds on their solutions. In the absence of the binary vector $y$ in (AP)($\xi$), training problems resemble those from superquantile-based supervised learning \cite{LaguelPillutlaMalickHarchaoui.21} and self-supervised learning \cite{LiuEtAl.21}. Our focus, however, is on settings {\em with} binary variables and this  lands us in a challenging territory. The training problems become large-scale mixed-binary problems, potentially nonlinear due to the choice of risk measure or organically because of $\psi_0, \psi_1, \dots, \psi_q$. The proposed decomposition algorithm leverages the (typical) near-separability of a training problem along the training set. A relaxation then leads to a fully-separable problem furnishing lower bounds on the minimum value of the training problem. This can be combined with any number of feasibility heuristics to construct upper bounds; we provide details for one approach. Preliminary numerical tests on a nonlinear binary optimization problem taken from target-search theory demonstrate that the training problems are solvable by the decomposition algorithm to reasonable accuracy in minutes.

Our focus on a finite training set $\{\xi^1, \dots, \xi^s\in \reals^r\}$ draws high-level connections with \cite{HadjiyiannisGoulartKuhn.11} and the general literature on supervised learning. It also raises question about how to ensure feasibility for the prescribed decisions at parameter vectors outside the training set, which we handle using a strategy resembling finite adaptivity as in \cite{HanasusantoKuhnWiesemann.15}. The local nature of our decision rules relates to those of the recent papers \cite{BertsimasShternSturt.22}, where each observed data point is associated with a {\em local} affine decision rule in the setting of linear optimization problems. Connections also emerge with general stability analysis of optimization problems. That area of study focuses on determining whether a parameterized problem has a set of (near-)minimizers that vary (Lipschitz) continuously under changing parameters; see the monographs \cite{VaAn,primer} as well as efforts based on metric regularity and calmness \cite{IoffeOutrata.08,Penot.10}, tilt-stability \cite{EberhardWenczel.12,LewisZhang.13,DrusvyatskiyLewis.13}, full-stability \cite{MordukhovichRockafellarSarabi.13}, and the truncated Hausdorff distance \cite{AttouchWets.93b,Royset.20b}. However, these studies are not tailored to binary decision variables and also lack a focus on decision rules. 

With the impressive successes of supervised learning over the last decades, one might wonder whether the following approach would suffice: (a) Generate $\{\xi^1, \dots, \xi^s\in \reals^r\}$. (b) For each $i=1, \dots, s$, obtain a minimizer $(x^i,y^i)$ of (AP)$(\xi^i)$. (c) Train a neural network on the data $\{\xi^i, (x^i, y^i), i=1, \dots, s\}$ to predict a minimizer of (AP)$(\xi)$ from an input vector $\xi$. Presently, this approach does not seem viable despite some successes in solving combinatorial optimization problems by neural networks \cite{KotaryFiorettoVanHentenryckWilder.21,LiuEtAl.21,BengioLodiProuvost.21,ShiZhang.22,Gasseetal.22} but this might change over time. Even if the approach can be carried out, it will be unable to provide the insight that can come from identifying an affine decision rule or similar ones with simple structure. Moreover, the benefit from a herculean effort to determine the mapping from $\xi$ to a minimizer of (AP)$(\xi)$ would anyhow be diminished by the fact that (AP)$(\xi)$ is just a {\em model} of a real-world decision process. Our experience is that a decision rule that prescribes good but potentially suboptimal decisions will be practically just as meaningful. In this paper, we develop such decision rules and examine their performance. Supervised learning of solution mappings is beyond the scope of the paper.  

Decision rules have traditionally been motivated by two- or multi-stage optimization problems under uncertainty (see, e.g., \cite{DelageIancu.15,GeorghiouKuhnWiesemann.19,YanikogluGorissenHertog.19}). The generality of (AP)($\xi$) makes the following results applicable to such settings in principle, but we omit a tailored development. Still, the numerical tests address search for moving targets across eight time periods. 

Our analysis is inherently local: Local consistency (Theorem \ref{thm:convergence}) ensures that decision rules produced by training points converging to $\bar \xi$ indeed prescribe a minimizer for (AP)($\bar\xi$) in the limit. Our nonasymptotic analysis (Theorem \ref{thm:error}) results in a suboptimality bound that is proportional to the size of any set of $\xi$-vectors under consideration. Thus, the results are especially well suited for perturbation (sensitivity) analysis involving ``small'' changes to $\xi$ and may require partitioning of any ``large'' set of $\xi$-vectors into smaller subsets (cf. \cite{VayanosKuhnRustem.11,PostekHertog.16,BertsimasDunning.16} for such approaches). We note, however, that the numerical tests show that the resulting decision rules tend to perform well across a wider range of $\xi$-vectors than predicted by our theory. 

The paper is organized as follows. Section 2 lays out the framework, defines training problems, and discusses how decision rules, if selected from a flexible enough class of mappings, recover minimizers of (AP)($\xi$). Section 3 establishes local consistency. Section 4 turns to a nonasymptotic analysis of decision rules for (AP)($\xi$) without the constraints $\psi_k(\xi,x,y) \leq 0$, $k=1, \dots, q$. Section 5 presents a decomposition algorithm for solving a large class of training problems and reports on preliminary numerical results. Section 6 gives conclusions. Section 7, intended for an online supplement, furnishes proofs, implementation details, and additional numerical results.

\section{Framework}\label{sec:framework}

We seek to construct a decision rule $(F,G) \in \cF \times \cG$ such that $(F(\xi),\bbH(G(\xi)))$ is a good decision for (AP)($\xi$), in some sense, for $\xi$ in a set $\Xi_0$. Traditionally, a focus has been on the {\em worst-case} values of $\psi_k(\xi,F(\xi),\bbH(G(\xi)))$, $k=0, 1, \dots, q$, across $\xi\in \Xi_0$ \cite{DelageIancu.15,GeorghiouKuhnWiesemann.19,YanikogluGorissenHertog.19}. Alternatively, one could consider the {\em average} performance of $(F,G)$ on a finite set $\Xi_0$, which is more aligned with statistical learning. We consider these two possibilities and many others within a unified framework by leveraging risk measures. Risk measures are widely used in all areas of decision making under uncertainty \cite{Royset.22b}. Averages and worst-case considerations then emerge as special cases, with the possibility of adapting risk measures to the particular needs of an application. Superquantiles are especially versatile tools to construct a variety of risk measures; see \cite{LaguelPillutlaMalickHarchaoui.21} for their use in machine learning.

\subsection{Training Problem}

We formalize the framework as follows. Given a finite set of {\em training points} $\xi^1, \dots, \xi^s\in \reals^r$, we seek to determine a decision rule that is ``good'' across these points. (The hope is that the decision rule is acceptable even beyond these training points, which indeed can be the case as we see below.) To connect with the theory of risk measures (see for example \cite{Royset.22b}) and facilitate analysis, we view the training points as the outcomes of a random vector $\bfxi$ defined on a {\em finite probability space} $(\Omega, 2^\Omega, P)$, where $P(\omega)>0$ for all $\omega\in \Omega$. That is, $\bfxi:\Omega\to \reals^r$ and $\{\bfxi(\omega)~|~\omega \in \Omega\} = \{\xi^k, \,k=1, \dots, s\}$. Throughout, we use boldface letters for random vectors.

In light of the actual problem, we seek $(F,G)\in \cF\times\cG$ such that $\psi_0(\bfxi, F(\bfxi), \bbH(G(\bfxi)))$ is low and each $\psi_k(\bfxi, F(\bfxi), \bbH(G(\bfxi)))$, $k=1, \dots, q$, is nonpositive, but these quantities are random variables making the meaning of ``low'' and ``nonpositive'' ambiguous. Risk measures allow us to define and adapt the meaning of such statements, while retaining computational tractability.

Following \cite[Definition 2.7]{Royset.22b}, a {\em measure of risk} (also called risk measure) $\cR$ assigns to a random variable $\bfeta$ a number $\cR(\bfeta) \in [-\infty,\infty]$ as a quantification of its {\em risk}, with this number being $\eta$ if $\prob\{\bfeta = \eta\} = 1$. Additional desirable properties associated with risk measures emerge below. 

For random variables $\bfeta$ (defined on $(\Omega, 2^\Omega, P)$), examples of risk measures include
\begin{align*}
  \cR(\bfeta) & = \Ex[\bfeta] = \sum_{\omega\in \Omega} P(\omega) \bfeta(\omega) &~~~\mbox{ (expectation)}\\
  \cR(\bfeta) & = \max_{\omega\in \Omega} \bfeta(\omega) &~~~\mbox{ (worst-case risk)}\\
  \cR(\bfeta) & = Q_\alpha(\bfeta) + \frac{1}{1-\alpha}\sum_{\omega\in \Omega} P(\omega) \max\big\{0, \bfeta(\omega) - Q_\alpha(\bfeta)\big\} &~~~\mbox{ (superquantile risk)}
\end{align*}
where $\alpha \in (0,1)$ and $Q_\alpha(\bfeta)$ is the $\alpha$-quantile of $\bfeta$; see, e.g., \cite{Royset.22b}.

The goal of low $\psi_0(\bfxi, F(\bfxi), \bbH(G(\bfxi)))$ and nonpositive $\psi_k(\bfxi, F(\bfxi), \bbH(G(\bfxi)))$, $k=1, \dots, q$, can now be interpreted as having low and nonpositive risks, respectively. For the risk measures $\cR_0, \cR_1, \dots, \cR_q$, this approach leads to the {\em training problem}
\begin{align}
\mbox{(TP)}~~~~~~\nnmin_{F\in \cF, G\in \cG}  ~\cR_0\Big(\psi_0\big(\bfxi,F(\bfxi),\bbH\big(G(\bfxi)\big)\big)\Big) &&&\nonumber\\
 \mbox{subject to } ~\cR_k\Big(\psi_k\big(\bfxi,F(\bfxi),\bbH\big(G(\bfxi)\big)\big)\Big) & \leq 0 && k=1, \dots, q\label{eqn:trainingproblemCon}\\
 \Big(F\big(\bfxi(\omega)\big), \bbH\big(G\big(\bfxi(\omega)\big)\big)\Big) & \in C &&\forall \omega\in \Omega\label{eqn:trainingproblemXcon}\\
 G\big(\bfxi(\omega)\big) &\in D_\epsilon(\delta) &&\forall \omega\in \Omega,\label{eqn:trainingproblemYcon}
\end{align}
where, for $\epsilon \in (0,\infty)$ and $\delta \in (\epsilon, \infty]$, we let
\[
D_\epsilon(\delta) = \big([-\delta, -\epsilon] \cup [\epsilon, \delta] \big)^m.
\]
The objective function in (TP) incentivizes a decision rule $(F,G)$ that produces a random variable $\psi_0(\bfxi,F(\bfxi),\bbH(G(\bfxi)))$ with a low risk (in the sense specified by $\cR_0)$. A conservative analyst might adopt the worst-case risk measure, while another one might simply consider the expectation. The constraints \eqref{eqn:trainingproblemCon} ensure that each random variable $\psi_k(\bfxi,F(\bfxi),\bbH(G(\bfxi)))$, $k=1, \dots, q$, has nonpositive risk in the sense of $\cR_k$. If $\cR_k$ is a superquantile risk measure at level $\alpha_k \in (0,1)$, then \eqref{eqn:trainingproblemCon} is equivalent to the buffered probability of $\psi_k(\bfxi,F(\bfxi),\bbH(G(\bfxi)))> 0$ being bounded from above by $1-\alpha_k$. This in turn is a conservative approximation of the chance constraint: $\prob(\psi_k(\bfxi,F(\bfxi),\bbH(G(\bfxi))) > 0) \leq 1-\alpha_k$; see \cite{RockafellarRoyset.10}. The constraints \eqref{eqn:trainingproblemXcon} ensure that $(F,G)$ produces a decision in $C$ regardless of the outcome of $\bfxi$.

The constraints \eqref{eqn:trainingproblemYcon} may seem superfluous but they play important practical and theoretical roles. Let the component functions of $G$ be written as $g_1, \dots, g_m$, i.e., $G(\xi) = (g_1(\xi), \dots, g_m(\xi))$. A requirement $G(\xi) \in D_\epsilon(\delta)$ amounts to having $g_i(\xi) \in [-\delta, -\epsilon] \cup [\epsilon, \delta]$ for each $i=1, \dots, m$. If $g_i(\xi) \in [-\delta, -\epsilon]$, then $\bbH(g_i(\xi)) = 0$, which prescribes the decision $y_i = 0$ in (AP)$(\xi)$. If $g_i(\xi)\in [\epsilon, \delta]$, then $\bbH(g_i(\xi)) = 1$, which prescribes the decision $y_i = 1$. In the absence of \eqref{eqn:trainingproblemYcon}, (TP) might produce a mapping $G$ with $g_i(\bfxi(\omega))$ being zero or arbitrarily close to zero for some $\omega\in \Omega$. This in turn would produce a ``fragile'' decision rule that prescribes a different decision for $y_i$ under $\bfxi(\omega)$ than for a nearly identical parameter vector $\xi$, i.e., $\bbH(g_i(\xi)) \neq \bbH(g_i(\bfxi(\omega)))$. This is a familiar situation from statistical learning. One would like the ``decision boundary'' in a classification problem to be some distance away from the location of the training data; $L_2$-regularized support vector machines make this requirement explicit by encouraging a wide ``margin;" see, e.g., \cite[Section 2.H]{primer}. The constraints \eqref{eqn:trainingproblemYcon} ensure that each $g_i$ has value no closer to zero than $\epsilon$ for all $\bfxi(\omega), \omega\in \Omega$. This produces a decision rule with a margin between the decision boundary and the training points. A larger $\epsilon$ tends to produce a wider margin, but its size is also affected by the flexibility of the mappings in $G$. The additional parameter $\delta$ could be infinity and thus having no effect in \eqref{eqn:trainingproblemYcon}. The restriction to a finite $\delta$ tends to be computationally beneficial (see Subsection \ref{subsec:implement}) and we include it as a possibility.

Given $\xi$ and a solution $(F,G)$ of (TP), the quality of $(F(\xi),\bbH(G(\xi)))$ compared to a minimizer of (AP)$(\xi)$ depends on several factors. In the following, we furnish theoretical and empirical results establishing that a near-minimizer of (TP) prescribes a nearly optimal decision for the actual problem (AP)$(\xi)$ as long as the training points are {\em near} $\xi$.

In implementation of a decision rule $(F,G)$, there might be a question about whether it is suitable for a particular parameter vector $\xi$. Is it close enough to the training points for the theoretical results to apply? While we provide a quantitative analysis in Section \ref{sec:error}, it is immediately clear that one can check whether $g_i(\xi) \in (-\epsilon,\epsilon)$ for some $i$, which would indicate that $\xi$ is close to the decision boundary. This would be of concern and might prompt abandoning the decision rule and solving (AP)$(\xi)$ directly. Likewise, violation of the constraint $(F(\xi),\bbH(G(\xi)))\in C$ can be checked. 

In practice, mappings in $\cF$ and $\cG$ are rarely flexible enough to track minimizers of (AP)$(\xi)$ as $\xi$ varies. However, if these collections are expansive enough, then we obtain the following relation between minimizers of the actual problem (AP)$(\xi)$ and those of the training problem (TP).

We recall that a risk measure $\cR$ is {\em monotone} if $\cR(\bfeta) \leq \cR(\bfeta')$ whenever $\bfeta(\omega) \leq \bfeta'(\omega)$ for every $\omega\in \Omega$. For example, expectation, worst-case, quantile, and superquantile risk measures are monotone \cite{Royset.22b}.

\begin{proposition}{\rm (solution recovery).}\label{prop:recovery} Suppose that each $\cR_1, \dots, \cR_q$ is the worst-case risk measure and there exists $(\bar F, \bar G) \in \cF \times \cG$ such that for every $\omega\in \Omega$ one has
\begin{equation*}\label{eqn:reconverycond}
\bar G\big(\bfxi(\omega)\big) \in D_\epsilon(\delta)~ \mbox{ and }~ \Big(\bar F\big(\bfxi(\omega)\big), \,\bbH\big(\bar G\big(\bfxi(\omega)\big)\big)\Big) \mbox{ is a minimizer of $(\rm{AP})(\bfxi(\omega))$.}
\end{equation*}
Then, the following hold:
\begin{enumerate}[(a)]

\item If the risk measure $\cR_0$ is monotone, then $(\bar F, \bar G)$ is a minimizer of $(\rm{TP})$.

\item If $\cR_0$ is the expectation and $(F^\star,G^\star)$ is a minimizer of $(\rm{TP})$, then for each $\omega\in \Omega$
\[
\Big(F^\star\big(\bfxi(\omega)\big), \,\bbH\big(G^\star\big(\bfxi(\omega)\big)\big)\Big) \mbox{ is a minimizer of $(\rm{AP})(\bfxi(\omega))$.}
\]
\end{enumerate}
\end{proposition}
\state Proof. See Supplementary Material in Section \ref{sec:supp}.\eop

The first assertion in the proposition establishes an optimality condition for the collection of problems  produced by setting $\xi = \bfxi(\omega)$, $\omega\in \Omega$, in (AP)$(\xi)$. It is necessary for minimizers of these problems (that satisfy the $D_\epsilon(\delta)$ condition) to also produce a minimizer of (TP). This holds for any monotone risk measure $\cR_0$. The second assertion in the proposition cannot be extended beyond the expectation assumption because (TP) tends to obscure or even eliminate the importance of (some) individual outcomes of $\bfxi$ when it involves a general risk measure in the objective function. In either assertion, $\cF\times\cG$ is assumed to be sufficiently rich to include mappings that can track minimizers of the actual problem (AP)$(\xi)$ as $\xi$ varies. In summary, there are circumstances when a decision rule obtained by solving a training problem prescribes minimizers of the actual problem for all values of $\xi$ corresponding to training points. Conversely, minimizers of the actual problem for $\xi$ across the training points may match a decision rule that solves the training problem. For neither part (a) nor part (b) are we able to consider a broader range of risk measures $\cR_1, \dots, \cR_q$ because a shift away from worst-case risk measures would break the equivalence between \eqref{eqn:trainingproblemCon} in (TP) and $\psi_k(\bfxi(\omega),F(\bfxi(\omega)),\bbH(G(\bfxi(\omega))))\leq 0$, $k=1, \dots, q$ in (AP)$(\bfxi(\omega))$ for all $\omega\in \Omega$.

While Proposition \ref{prop:recovery} connects fundamentally the training problem (TP) to the actual problem, we need to go further and examine their relation when $\cF\times\cG$ is less expansive, for example, containing only affine decision rules. We also would like to consider a wider array of risk measures as well as approximations in (TP), possibly computationally motivated. Sections \ref{sec:conv} and \ref{sec:error} address these concerns.

\subsection{Implementation and Solution of Training Problems}\label{subsec:implement}

The choice of a computational approach for solving (TP) depends on the representation of mappings in $\cF$ and $\cG$, the properties of $\psi_0, \psi_1, \dots, \psi_q$, the nature of $C$, the risk measures, and other factors. The Heaviside function $\bbH$ causes particular challenges, but recent advances in computational nonsmooth optimization such as in \cite{PenaLuedtkeWachter.20,CuiLiuPang.22} might lead to viable algorithms for (TP). In Section \ref{sec:num}, we rely on integer programming solvers and then the following trivial fact is useful: With $g_i$ being the $i$th component function of $G\in \cG$, $\delta \in (\epsilon, \infty)$, and $\xi\in \reals^r$, one has
\[
\bbH\big(G(\xi)\big) = y, ~G(\xi) \in D_\epsilon(\delta)  ~~~\Longleftrightarrow ~~-\delta + (\delta+\epsilon)y_i \leq g_i(\xi) \leq -\epsilon + (\delta + \epsilon)y_i, ~ y_i \in \{0,1\}, ~ \forall i.
\]
This fact allows for reformulations of (TP). For affine decision rules, optimizing over $(F,G) \in \cF\times \cG$ is equivalent to optimizing over $A\in \reals^{n\times r}$, $B\in \reals^{m\times r}$, $a\in\reals^n$, $b\in \reals^m$, with these matrices and vectors being (potentially) restricted to some set $S$. As an example, suppose that $\cR_0, \cR_1, \dots, \cR_q$ are superquantile measures of risk at level $\alpha_0, \alpha_1, \dots, \alpha_q\in (0,1)$, respectively. Then, leveraging the Rockafellar-Uryasev formula for superquantiles (cf. \cite{RockafellarUryasev.13} or \cite{Royset.22b} for details), (TP) can be reformulated as:
\begin{align*}
\mbox{(TP-super)}~~\nnmin_{A,a,B,b,u,y,z}  ~z_0 + \frac{1}{1-\alpha_0} \sum_{\omega\in \Omega} P(\omega) u_0(\omega)& & &\\
 \mbox{subject to } ~ z_k + \frac{1}{1-\alpha_k} \sum_{\omega\in \Omega} P(\omega) u_k(\omega) & \leq 0 &&k=1, \dots, q\\
 \psi_k\big(\bfxi(\omega),A\bfxi(\omega) + a,y(\omega)\big) - z_k & \leq u_k(\omega) &&\forall \omega\in \Omega, ~k=0, 1,  \dots, q\\
                           \big(A\bfxi(\omega) + a, y(\omega)\big) & \in C &&\forall \omega\in \Omega\\
                           -\delta + (\delta+\epsilon)y_i(\omega) \leq \big\langle B_i, \bfxi(\omega)\big\rangle + b_i & \leq -\epsilon + (\delta + \epsilon)y_i(\omega) && \forall \omega\in \Omega, ~i=1, \dots, m\\
                           (A,a,B,b) \in S,~~~~ & (z_0, z_1, \dots, z_q) \in \reals^{1+q}&&\\
                           \big(u_0(\omega), u_1(\omega), \dots, u_q(\omega)\big) \in [0,\infty)^{1+q}, ~ & \big(y_1(\omega), \dots, y_m(\omega)\big)  \in \{0,1\}^{m} && \forall \omega\in \Omega.
\end{align*}
Here, $B_i$ is the $i$th row of the matrix $B$ and $b_i$ is the $i$th element of the vector $b$. If the set $S\subset \reals^{n\times r} \times \reals^n \times \reals^{m\times r} \times \reals^m$ is convex and each $\psi_k(\bfxi(\omega),\cdot,\cdot)$, $k=0, 1, \dots, q$, is a convex function for all $\omega\in \Omega$, then this reformulation has a convex continuous relaxation provided that $C=C' \cap (\reals^n \times \{0,1\}^m)$ for some convex set $C'\subset\reals^{n+m}$. Even if well-structured, the reformulation tends to be large-scaled and this is a common characteristic of (TP). We develop a decomposition algorithm in Section \ref{sec:num} that solves practically important instances.

\section{Local Consistency}\label{sec:conv}

While the purpose of a decision rule $(F,G)$ is to prescribe a decision $(F(\xi),\bbH(G(\xi)))$ for (AP)$(\xi)$ at a collection of future, unknown parameter vectors $\xi$, we still would like a local consistency guarantee at a particular $\bar \xi$ of the following form: If the training points $\xi^1, \dots, \xi^s$ tend to $\bar\xi$, then the corresponding decision rules produced by (TP) can only converge to a decision rule $(F^\star,G^\star)$ for which $(F^\star(\bar\xi),\bbH(G^\star(\bar\xi)))$ is a minimizer of the actual problem (AP)$(\bar\xi)$. Informally, this means that solving (TP) using a finite set of training points near $\bar \xi$ produces a decision rule $(F,G)$ with a near-optimal prescription at $\bar \xi$. In this section, we present three results formalizing this claim: Theorem \ref{thm:convergence} states the main result under the assumption of monotone risk measures $\cR_1, \dots, \cR_q$. Corollary \ref{cor:convergenceNoIneq} specializes the theorem to cases without inequality constraints in (AP)($\xi$), and this eliminates all assumptions about $\cR_1, \dots, \cR_q$. Corollary \ref{cor:convergence} returns to the setting with inequality constraints and reports a one-sided bound without the need for monotonicity assumptions on $\cR_1, \dots, \cR_q$.

To formally address a situation with changing training points, we again consider a finite probability space $(\Omega,2^\Omega,P)$, with $P(\omega)>0$ for all $\omega\in \Omega$. Let $\nats = \{1, 2, \dots\}$ be the natural numbers and $\{\bar\bfxi,\bfxi^\nu:\Omega\to \reals^r, \nu\in\nats\}$ be random vectors defined on the probability space. (Throughout, $\nu$ indexes sequences.) The random vector $\bar\bfxi$ has the constant value $\bar\xi$, i.e., $\bar\bfxi(\omega) = \bar \xi$ for all $\omega\in \Omega$. The random vector $\bfxi^\nu$ represents a set of training points, with the outcomes $\bfxi^\nu(\omega), \omega\in \Omega$ being those points. Thus, pointwise convergence of the random vectors $\bfxi^\nu$ to $\bar\bfxi$, written as $\bfxi^\nu\to \bar\bfxi$, models situations with training points approaching $\bar\xi$. (We note that since the probability space is finite and $P(\omega)>0$ for all $\omega \in \Omega$, pointwise convergence is equivalent to almost-sure convergence and also $\cL^2$-convergence.)

As compared to (TP), we introduce three refinements that address adjustments possibly needed in practical implementation. First, we include {\em regularization} as represented by a function $h^\nu:\cF\times\cG\to [0,\infty)$, which might change with $\nu$ or simply be the zero function. We assume throughout that these functions are nonnegatively valued. While the choices of $\cF$ and $\cG$ define admissible mappings, regularization can be used to encourage particular characteristics such as sparsity of affine mappings. Second, the functions $\psi_0, \psi_1, \dots, \psi_q$ might be unavailable or undesirable and we permit their replacement by {\em approximating functions} $\psi_k^\nu:\reals^r\times\reals^n\times\reals^m\to \reals$, $k= 0, 1, \dots, q$, which we assume throughout to be continuous. For example, $\psi_k^\nu$ might be a more desirable, smooth approximation of $\psi_k$, or it could represent the fact that $\psi_k$ can only be evaluated imprecisely. Third, we allow $\delta$ in \eqref{eqn:trainingproblemYcon} to vary with $\nu$. For example, $\delta$ could be $\infty$ but implementations in (TP-super) require finite values. Thus, one might consider a sequence $\{\delta^\nu, \nu\in\nats\}$ of such parameters tending to $\delta = \infty$.

These adjustments lead to a {\em sequence of training problems}: For each $\nu\in\nats$, we define
\begin{align}
\mbox{(TP)}^\nu~~~~~~\nnmin_{F\in \cF, G\in \cG}  ~\cR_0\Big(\psi_0^\nu\big(\bfxi^\nu,F(\bfxi^\nu),\bbH\big(G(\bfxi^\nu)\big)\big)\Big) & + h^\nu(F,G) &&\nonumber\\
 \mbox{subject to } ~~\cR_k\Big(\psi_k^\nu\big(\bfxi^\nu,F(\bfxi^\nu),\bbH\big(G(\bfxi^\nu)\big)\big)\Big) & \leq 0, && k=1, \dots, q\label{eqn:riskconstr}\\
                           \Big(F\big(\bfxi^\nu(\omega)\big), \bbH\big(G\big(\bfxi^\nu(\omega)\big)\big)\Big) & \in C && \forall \omega\in \Omega\nonumber\\
                           G\big(\bfxi^\nu(\omega)\big) &\in D_\epsilon(\delta^\nu) && \forall \omega\in \Omega.\label{eqn:Dcon}
\end{align}

To study the convergence of sequences of decision rules, we equip $\cF\times\cG$ with the topology of uniform convergence. Let $\Xi\subset \reals^r$ be a compact set that contains $\bar\xi$ as well as $\{\bfxi^\nu(\omega), \nu\in\nats, \omega\in \Omega\}$. For the remainder of the paper, we suppose that
\[
  \cF \subset \{F:\Xi\to \reals^n~|~F \mbox{ continuous } \} ~~\mbox{ and } ~~\cG \subset \{G:\Xi\to \reals^m~|~G \mbox{ continuous } \}.
\]
The restriction to continuous decision rules is hardly significant in practice; affine, piecewise affine, and polynomial decision rules are continuous and the same holds for common neural nets. For mappings in $\cF\times \cG$, we write $(F^\nu,G^\nu)\to (F,G)$ when $(F^\nu,G^\nu)$ converges uniformly to $(F,G)$. (Such convergence is equivalent to having $\max\{\max_{\xi\in \Xi} \|F^\nu(\xi)-F(\xi)\|_\infty, \max_{\xi\in \Xi} \|G^\nu(\xi)-G(\xi)\|_\infty\}\to 0$.) When $\cF$ and $\cG$ consist of affine mappings and $\Xi\subset\reals^r$ is full dimensional, we note that uniform convergence takes place if and only if the coefficients that define the affine mappings converge.

Regardless of the type of convergence, if $\alpha^\nu$ converges to $\alpha$ (merely) along a subsequence of indices $N\subset\nats$, then we write $\alpha^\nu \Nto \alpha$.

We also need the following terminology. For $\tau^\nu \in [0,\infty)$, we say that the decision rule $(F^\nu,G^\nu)\in \cF\times\cG$ is a $\tau^\nu$-{\em minimizer} of $\mbox{(TP)}^\nu$ if it satisfies the constraints in that problem and
\[
\cR_0\Big(\psi_0^\nu\big(\bfxi^\nu,F^\nu(\bfxi^\nu),\bbH\big(G^\nu(\bfxi^\nu)\big)\big)\Big) + h^\nu(F^\nu,G^\nu) \leq \mathfrak{m}^\nu + \tau^\nu < \infty,
\]
where $\mathfrak{m}^\nu$ is the minimum value of  $\mbox{(TP)}^\nu$. Similarly, with $\mathfrak{m}(\xi)$ being the minimum value of $\mbox{(AP)}(\xi)$ and $\tau \in [0,\infty)$, we define
\begin{align*}
&\tau\mbox{-}\nargmin_{(x,y)\in C} \big\{\psi_0(\xi,x,y)~\big|~\psi_k(\xi,x,y)\leq 0, k=1, \dots, q\big\}\\
& = \big\{(x,y)\in C~|~\psi_0(\xi,x,y) \leq \mathfrak{m}(\xi) +\tau < \infty,~ \psi_k(\xi,x,y)\leq 0, k=1, \dots, q\big\}.
\end{align*}

We denote by $\nt S$ the {\em interior} of a set $S\subset\reals^m$. We recall that a risk measure $\cR$ is {\em continuous} when $\bfeta^\nu\to \bfeta$ implies $\cR(\bfeta^\nu)\to \cR(\bfeta)$. The (Euclidean) {\em balls} on $\reals^n$ are $\ball(x,\rho) = \{x' \in \reals^n~|~\|x'-x\|_2 \leq \rho\}$. 

We are then ready to state the main theorem of this section, which establishes local consistency of decision rules in quite general settings. It is followed by two corollaries dealing with more specific situations.

\begin{theorem}{\rm (local consistency).}\label{thm:convergence} Suppose that $\delta^\nu \in (\epsilon,\infty] \to \delta \in (\epsilon, \infty]$, $\tau^\nu\in [0,\infty)\to 0$, $\bfxi^\nu\to \bar\bfxi$, and the following hold:
\begin{enumerate}[(a)]

\item The risk measures $\cR_0, \cR_1, \dots, \cR_q$ are real-valued and continuous; $\cR_1, \dots, \cR_q$ are also monotone.

\item For any $k=0, 1, \dots, q$, $\xi\in \Xi$, and $(x,y)\in C$, one has $\psi_k^\nu(\xi^\nu,x^\nu,y)\to\psi_k(\xi,x,y)$ whenever $\xi^\nu\in \Xi\to \xi$ and $(x^\nu,y)\in C\to (x,y)$.

\item $F+a\in\cF$ for any $F\in \cF$ and $a\in \reals^n$; and $G+b\in\cG$ for any $G\in \cG$ and $b\in \reals^m$.

\item There exists $h:\cF\times\cG\to [0,\infty)$ with the properties that $h(F,G) = h(F+a,G+b)$ for any $F\in \cF$, $a\in \reals^n$, $G\in \cG$, and $b\in \reals^m$; and $h^\nu(F^\nu,G^\nu)\to h(F,G)$ whenever $(F^\nu,G^\nu)\to (F,G)$.

\item There exists $(F,G)\in \cF\times\cG$ producing $G(\bar\xi) \in D_\epsilon(\delta)$ and a feasible point $(F(\bar\xi),\bbH(G(\bar\xi)))$ for $(\rm{AP})(\bar\xi)$. For each $(F,G)$ with these properties there exist $x^\nu\to F(\bar \xi)$ and $\gamma^\nu \in (0,\infty)\to 0$ that also satisfy
\begin{align*}
\psi_k\Big(\bar\xi,x^\nu,\bbH\big(G(\bar\xi)\big)\Big) &< 0, ~k=1, \dots, q, ~~~\nu\in\nats\\
\Big(z^\nu, \bbH\big(G(\bar\xi)\big)\Big) &\in C ~~\forall z^\nu \in \ball(x^\nu,\gamma^\nu), ~~~\nu\in\nats.
\end{align*}

\item For each $\nu\in \nats$, $(\hat F^\nu,\hat G^\nu)$ is a $\tau^\nu$-minimizers of $(\rm{TP})^\nu$.

\end{enumerate}
If there is a subsequence $N\subset\nats$ such that $(\hat F^\nu, \hat G^\nu)\Nto (F^\star,G^\star)$, then $\bbH(\hat G^\nu(\bfxi^\nu(\omega))) = \bbH(G^\star(\bar \xi))$ for all $\omega\in \Omega$ when $\nu\in N$ is sufficiently large, $G^\star(\bar \xi) \in D_\epsilon(\delta)$,
\[
\Big(F^\star(\bar \xi), \,\bbH\big(G^\star(\bar \xi)\big)\Big) \in \nargmin_{(x,y)\in C} \big\{\psi_0(\bar \xi,x,y) \,\big|\, \psi_k(\bar \xi,x,y) \leq 0, ~k=1, \dots, q\big\},
\]
and
\begin{align*}
& \cR_0\Big(\psi_0^\nu\big(\bfxi^\nu,\hat F^\nu(\bfxi^\nu),\bbH\big(\hat G^\nu(\bfxi^\nu)\big)\big)\Big) + h^\nu(\hat F^\nu,\hat G^\nu)\\
& \Nto  \inf_{(x,y)\in C} \big\{\psi_0(\bar \xi,x,y) ~\big| ~\psi_k(\bar \xi,x,y) \leq 0, k=1, \dots, q\big\} + h(F^\star,G^\star).
\end{align*}
Under the additional assumption that $\cF$ and $\cG$ contain their respective zero mappings and $h(F,G)=0$ only when $F$ and $G$ are constant mappings, we also have that $F^\star$ and $G^\star$ are constant mappings and, consequently, $h(F^\star,G^\star)=0$.
\end{theorem}
\state Proof. See Supplementary Material in Section \ref{sec:supp}.\eop

The theorem furnishes sufficient conditions for decision rules generated by (TP)$^\nu$ to provide in the limit an optimal decision for the actual problem (AP)$(\bar\xi)$ given a specific $\bar\xi\in\reals^r$. Thus, it gives the fundamental assurance that (TP)$^\nu$ cannot produce a decision rule that is ``poor'' for (AP)$(\bar \xi)$ if it uses training points ``near'' $\bar \xi$. In fact, the theorem establishes that the ``correct'' $y$-decision is achieved using training points in a neighborhood of $\bar \xi$, i.e., for some finite $\nu$. Additionally, the objective function values from (TP)$^\nu$ provide estimates of the minimum value of (AP)($\bar\xi$). We note that the number of training points remains constant in the theorem, but the points approach $\bar\xi$ as $\nu\to \infty$.

The conditions of Theorem \ref{thm:convergence} are generally mild. Assumption (a) is commonly satisfied such as in the case of expectation, quantiles, superquantiles, and worst-case risk measures. We may also have $\cR_0(\bfeta) = \Ex[\bfeta] + \lambda\var(\bfeta)$ or $\cR_0(\bfeta) = \Ex[\bfeta] + \lambda\std(\bfeta)$ for any positive $\lambda$, where $\var(\bfeta)$ and $\std(\bfeta)$ are the variance and standard deviation of the random variable $\bfeta$, respectively, but such risk measures cannot serve in the constraints as they are not monotone. Assumption (b) holds automatically if there is no approximation (i.e., $\psi_k^\nu = \psi_k$) because $\psi_k$ is continuous. In particular, $\psi_k^\nu$ and $\psi_k$ can be nonconvex and nonsmooth. Assumption (c) states that mappings in $\cF$ and $\cG$ can be shifted up or down by a constant vector and remain admissible, which is not restrictive in practice. Assumption (d) might be void if there is no regularization (TP)$^\nu$. But, it is also satisfied for affine mappings of the form $F(\xi) = A\xi + a$ and $G(\xi) = B\xi + b$ when $\Xi\subset\reals^r$ is full dimensional and
\begin{equation}\label{eqn:hexample}
h^\nu(F,G) = h(F,G) = \theta_1 \|A\| + \theta_2 \|B\|,
\end{equation}
where $\theta_1,\theta_2 \in [0,\infty)$ and the matrix norms are arbitrary. This produces a reasonable regularization term in (TP)$^\nu$ that encourages slowly varying mappings. Assumption (d) holds in this case because, for $F^\nu(\xi) = A^\nu\xi + a^\nu$ and $G^\nu(\xi) = B^\nu\xi + b^\nu$, one has 
\[
(F^\nu,G^\nu) \to (F,G) ~~\Longrightarrow~~h^\nu(F^\nu,G^\nu) = \theta_1 \|A^\nu\| + \theta_2 \|B^\nu\| \to \theta_1 \|A\| + \theta_2 \|B\| = h(F,G).
\]
Assumption (e) ensures that $\cF\times\cG$ can produce a feasible decision in the actual problem (AP)$(\bar\xi)$. In the absence of the constraints $\psi_k(\xi,x,y)\leq 0$, $k=1, \dots, q$, in (AP)$(\xi)$ and therefore also no risk constraints \eqref{eqn:riskconstr}, assumption (e) simply amounts to this feasibility requirement, now pertaining only to $C$. (Corollary \ref{cor:convergenceNoIneq} simplifies this even further.) If there are constraints $\psi_k(\xi,x,y)\leq 0$, $k=1, \dots, q$, we also assume a Slater-type condition (cf. \cite[Example 5.47]{primer}) for the actual problem: every feasible point in (AP)$(\bar\xi)$ can be approached from ``inside'' in a sense made specific in assumption (e). Without this property, there is no clear path for (TP)$^\nu$ to produce decision rules that are feasible in the actual problem, even in the limit. As specified by assumption (f), the training problem (TP)$^\nu$ can be solved with a nonzero tolerance.

The additional result at the end of Theorem \ref{thm:convergence} imposes the assumption that $\cF$ and $\cG$ contain the zero functions, i.e., mappings $F,G$ with $F(\xi) = 0$ and $G(\xi) = 0$ for all $\xi\in \Xi$. This is typically not restrictive. The assumption that $h(F,G) = 0$ only when $F$ and $G$ are constant mappings holds, for example, if $h$ is as in \eqref{eqn:hexample}.

For the setting without inequality constraints in (AP)($\xi$) and without approximations of cost or regularization terms, Theorem \ref{thm:convergence} simplifies to the following corollary. 

\begin{corollary}{\rm (local consistency without inequality constraints).}\label{cor:convergenceNoIneq} Suppose that $\delta^\nu \in (\epsilon,\infty] \to \delta \in (\epsilon, \infty]$, $\tau^\nu\in [0,\infty)\to 0$, $\bfxi^\nu\to \bar\bfxi$, and the following hold:
\begin{enumerate}[(a)]

\item The risk measure $\cR_0$ is real-valued and continuous.

\item $F+a\in\cF$ for any $F\in \cF$ and $a\in \reals^n$; and $G+b\in\cG$ for any $G\in \cG$ and $b\in \reals^m$.

\item The function $h:\cF\times\cG\to [0,\infty)$ is continuous with the property that $h(F,G) = h(F+a,G+b)$ for any $F\in \cF$, $a\in \reals^n$, $G\in \cG$, and $b\in \reals^m$.

\item For each $\nu\in \nats$, $(\hat F^\nu,\hat G^\nu)$ is a $\tau^\nu$-minimizers of $(\rm{TP})^\nu$ with $h^\nu = h$ and $\psi_0^\nu = \psi_0$, but without \eqref{eqn:riskconstr}.

\end{enumerate}
If there is a subsequence $N\subset\nats$ such that $(\hat F^\nu, \hat G^\nu)\Nto (F^\star,G^\star)$, then 
\[
\Big(F^\star(\bar \xi), \,\bbH\big(G^\star(\bar \xi)\big)\Big) \in \nargmin_{(x,y)\in C} \psi_0(\bar \xi,x,y).
\]
\end{corollary}
\state Proof. The result follows directly from Theorem \ref{thm:convergence} after recalling that $\psi_0$ is continuous by assumption and after constructing a feasible point for (AP)($\bar \xi$) to satisfy assumption (e) of that theorem. For the latter, let $\hat\omega\in \Omega$ be arbitrary and define 
\[
F(\xi) = \hat F^\nu(\xi) + \hat F^\nu\big(\bfxi^\nu(\hat\omega)\big) - \hat F^\nu(\bar\xi) ~~\mbox{ and } ~~G(\xi) = \hat G^\nu(\xi) + \hat G^\nu\big(\bfxi^\nu(\hat\omega)\big) - \hat G^\nu(\bar\xi) ~~\mbox{ for all } \xi\in \Xi.
\]
Then, $(F,G)\in(\cF,\cG)$ and satisfies $G(\bar \xi) \in D_\epsilon(\delta)$ as well as $(F(\bar \xi), \bbH(G(\bar \xi))) \in C$.\eop

In the absence of inequality constraints in the actual problem, the corollary establishes that essentially all meaningful risk measures (including expectation, quantiles, superquantiles, worst-case, mean-plus-standard-deviation, mean-plus-variance, and various combinations thereof) and nearly any families of continuous mappings lead to training problems with local consistency.

Even in the original setting of Theorem \ref{thm:convergence}, the assumptions can be relaxed considerably if we only seek an upper bound on the performance in the actual problem (AP)$(\bar\xi)$ of a decision rule generated by (TP)$^\nu$.

\begin{corollary}{\rm (limiting bound).}\label{cor:convergence} Suppose that $\delta^\nu \in (\epsilon,\infty] \to \delta \in (\epsilon, \infty]$, $\bfxi^\nu\to \bar\bfxi$, and the following hold:
\begin{enumerate}[(a)]

\item The risk measures $\cR_0, \cR_1, \dots, \cR_q$ are real-valued and continuous.

\item For any $k=0, 1, \dots, q$, $\xi\in \Xi$, and $(x,y)\in C$, one has $\psi_k^\nu(\xi^\nu,x^\nu,y)\to\psi_k(\xi,x,y)$ whenever $\xi^\nu\in \Xi\to \xi$ and $(x^\nu,y)\in C\to (x,y)$.

\item There exists $h:\cF\times\cG\to [0,\infty)$ with the property that $\nliminf h^\nu(F^\nu,G^\nu)\geq h(F,G)$ whenever $(F^\nu,G^\nu)\to (F,G)$.

\item For each $\nu\in \nats$, $(\hat F^\nu,\hat G^\nu)$ is feasible for $(\rm{TP})^\nu$ and
\[
\cR_0\Big(\psi_0^\nu\big(\bfxi^\nu,\hat F^\nu(\bfxi^\nu),\bbH\big(\hat G^\nu(\bfxi^\nu)\big)\big)\Big) + h^\nu(\hat F^\nu,\hat G^\nu) \leq \gamma^\nu.
\]

\end{enumerate}
If there is a subsequence $N\subset\nats$ such that $(\hat F^\nu, \hat G^\nu)\Nto (F^\star,G^\star)$ and $\nliminf_{\nu\in N} \gamma^\nu<\infty$, then
\[
\Big(F^\star(\bar\xi),\, \bbH\big(G^\star(\bar\xi)\big)\Big) ~~\mbox{ is feasible for $(\rm{AP})(\bar\xi)$, }
\]
$\bbH(\hat G^\nu(\bfxi^\nu(\omega)))$ $=$ $\bbH(G^\star(\bar \xi))$ for all $\omega\in \Omega$ when $\nu\in N$ is sufficiently large, $G^\star(\bar \xi) \in D_\epsilon(\delta)$, and
\[
\psi_0\Big(\bar \xi,F^\star(\bar\xi),\bbH\big(G^\star(\bar\xi)\big)\Big) + h(F^\star,G^\star) \leq \nliminf_{\nu\in N} \gamma^\nu.
\]
\end{corollary}
\state Proof. See Supplementary Material in Section \ref{sec:supp}.\eop

The corollary establishes that {\em any feasible} decision rule $(\hat F^\nu,\hat G^\nu)$ in the training problem (TP)$^\nu$, regardless of whether is it optimal, near-optimal, or neither, can be used to conservatively bound the performance of a limiting decision rule $(F^\star,G^\star)$ when implemented for the actual problem (AP)$(\bar\xi)$. In particular, the limiting decision rule is feasible for the actual problem  (including with any number of inequality constraints) provided that $\nliminf_{\nu\in N}\gamma^\nu<\infty$.  The values $\gamma^\nu$ produced by the training problems furnish in the limit an upper bound
\[
\psi_0\Big(\bar \xi,F^\star(\bar\xi),\bbH\big(G^\star(\bar\xi)\big)\Big) \leq \nliminf_{\nu\in N} \gamma^\nu
\]
because $h(F^\star,G^\star)\geq 0$. Of course, more optimization effort applied to (TP)$^\nu$ would tend to produce a better $(\hat F^\nu,\hat G^\nu)$ and thus a lower $\gamma^\nu$. This translates into a better bound $\nliminf_{\nu\in N}\gamma^\nu$ on how well the limiting decision rule $(F^\star,G^\star)$ performs on (AP)$(\bar\xi)$.

\section{Nonasymptotic Error Bounds}\label{sec:error}

We next turn to a quantification of the level of suboptimality for a decision rule produced by the training problem (TP)$^\nu$ given {\em  any fixed} $\nu\in\nats$, when applied to the actual problem (AP)$(\xi)$ with $\xi$ in some set $\Xi$. In contrast to Section \ref{sec:conv}, the analysis is nonasymptotic and establishes {\em uniform} performance guarantees across $\Xi$. Thus, we think of the index $\nu$ being fixed at some value. We follow the setting of Section \ref{sec:conv}, but now $\Xi\subset\reals^r$ is a compact set containing $\bfxi^\nu(\omega)$, $\omega\in \Omega$ for the specific $\nu$ under consideration. We provide two theorems. The first one furnishes a lower bound on the minimum values of (AP)$(\xi)$ for $\xi\in \Xi$ expressed in terms on the minimum value of (TP)$^\nu$. The second theorem provides near-optimality guarantees for three decision rules produced by (TP)$^\nu$.  In both cases, we assume monotone and convex risk measures.

Throughout this section, we omit the constraints $\psi_k(\xi,x,y)\leq 0$, $k=1, \dots, q$, from (AP)$(\xi)$ and the risk constraints \eqref{eqn:riskconstr} from (TP)$^\nu$. As indicated by assumption (e) in Theorem \ref{thm:convergence}, these constraints cause significant complications because for any fixed $\nu$ the training problem  (TP)$^\nu$ might be infeasible and/or provide little information about what constitute a feasible decision in (AP)$(\xi)$ for $\xi$ outside the set of training points. Thus, we assume in this section that constraints dependent on $\xi$ in (AP)$(\xi)$ have been incorporated into the objective function $\psi_0$ using penalty expressions or elastic programming; see \cite{BrownGraves.75,Royset.21} for strong practical and theoretical justifications for such model formulations. (Constraints independent of $\xi$ remain as captured by the set $C$.)

We adopt the notation
\[
\diam \, \Xi = \max_{\xi,\xi'\in \Xi} \|\xi - \xi'\|_2.
\]
Recall that a risk measure $\cR$ is {\em convex} if $\cR((1-\lambda)\bfeta + \lambda \bfeta') \leq (1-\lambda)\cR(\bfeta) + \lambda \cR(\bfeta')$ for every $\lambda\in [0,1]$ and random variables $\bfeta$, $\bfeta'$. Expectations, worst-case risk measures, and superquantile risk measures are convex,  but those defined by quantiles are generally not (see, e.g., \cite[Example 8.8]{primer} for a counterexample).

The first theorem is a building block for the second one, but also stands on its own by providing a lower bound.

\begin{theorem}{\rm (lower bound).}\label{thm:lower} For fixed $\nu\in\nats$, suppose that $\delta,\delta^\nu \in (\epsilon,\infty]$ and the following hold:
\begin{enumerate}[(a)]

\item The risk measure $\cR_0$ is real-valued, monotone, and convex.

\item $h^\nu(F,G) = 0$ when $F$ and $G$ are constant mappings.

\item There are $\sigma^\nu, \kappa_0\in [0,\infty)$ such that for all $(x,y)\in C$ one has
\begin{align*}
\big|\psi_0(\xi,x,y) - \psi_0(\xi',x,y)\big| & \leq \kappa_0 \|\xi - \xi'\|_2 ~~~\forall \xi,\xi'\in \Xi\\
\big|\psi_0^\nu(\xi,x,y) - \psi_0(\xi,x,y)\big| & \leq \sigma^\nu ~~ \forall \xi\in \Xi.
\end{align*}

\item $F+a\in\cF$ for any $F\in \cF$ and $a\in \reals^n$; $G+b\in\cG$ for any $G\in \cG$ and $b\in \reals^m$; and both $\cF$ and $\cG$ contain their respective zero functions.

\end{enumerate}
If \,$\mathfrak{m}^\nu$ is the minimum value of $(\rm{TP})^\nu$ without risk constraints \eqref{eqn:riskconstr}, then
\[
\mathfrak{m}^\nu  - \sigma^\nu - \kappa_0 \,\diam \, \Xi ~\leq \inf_{(x,y)\in C} \psi_0(\xi,x,y)~~~ \forall \xi\in \Xi.
\]
\end{theorem}
\state Proof. Let $\bar\xi\in \Xi$. If $\inf_{(x,y)\in C} \psi_0(\bar\xi,x,y) = \infty$, then the claimed inequality holds trivially for $\xi = \bar\xi$. Let $\gamma \in (0,\infty)$. Next, suppose that $\inf_{(x,y)\in C} \psi_0(\bar\xi,x,y) \in\reals$. There exists $(x(\bar\xi),y(\bar\xi))\in C$ such that
\[
\psi_0\big(\bar\xi,x(\bar\xi),y(\bar\xi)\big) \leq \inf_{(x,y)\in C} \psi_0(\bar\xi,x,y) + \gamma.
\]
Let $y_i(\bar\xi)$ be the $i$th component of $y(\bar\xi)$. We construct $F^\nu:\Xi\to \reals^n$ by setting $F^\nu(\xi) = x(\bar\xi)$ for all $\xi\in \Xi$ and $G^\nu:\Xi\to \reals^m$ by defining its $i$th component function as
\[
\forall \xi\in \Xi:~~~~~g_i^\nu(\xi) = \begin{cases}
\epsilon & \mbox{ if } y_i(\bar\xi)=1\\
-\epsilon & \mbox{ if } y_i(\bar\xi)=0.
\end{cases}
\]
Consequently, $F^\nu\in \cF$ and $G^\nu\in \cG$ by virtue of being constant mapping; see assumption (d). For any $\omega\in \Omega$, $G^\nu(\bfxi^\nu(\omega)) = G^\nu(\bar\xi) \in D_\epsilon(\delta^\nu)$ by construction. We also have that $F^\nu(\bfxi^\nu(\omega)) = x(\bar\xi)$ and $\bbH(G^\nu(\bfxi^\nu(\omega))) = y(\bar\xi)$. It trivially follows that $(F^\nu(\bfxi^\nu(\omega)), \bbH(G^\nu(\bfxi^\nu(\omega)))) \in C$. This means that $(F^\nu,G^\nu)$ is feasible for (TP)$^\nu$ without risk constraints \eqref{eqn:riskconstr} and
\begin{align*}
  \mathfrak{m}^\nu & \leq \cR_0\Big(\psi_0^\nu\big(\bfxi^\nu,F^\nu(\bfxi^\nu),\bbH\big(G^\nu(\bfxi^\nu)\big)\big)\Big) + h^\nu(F^\nu,G^\nu) - \psi_0\big(\bar\xi,x(\bar\xi),y(\bar\xi)\big) + \psi_0\big(\bar\xi,x(\bar\xi),y(\bar\xi)\big)\\
  & \leq \cR_0\Big(\psi_0^\nu\big(\bfxi^\nu,F^\nu(\bfxi^\nu),\bbH\big(G^\nu(\bfxi^\nu)\big)\big) - \psi_0\big(\bar\xi,x(\bar\xi),y(\bar\xi)\big)\Big)  + \inf_{(x,y)\in C} \psi_0(\bar\xi,x,y) + \gamma,
\end{align*}
where we use the facts that $h^\nu(F^\nu,G^\nu)=0$ because $(F^\nu,G^\nu)$ are constant mappings (cf. assumption (b)) and $\cR_0(\bfeta' + \bfeta) = \cR_0(\bfeta') + \eta$ for any random variables $\bfeta'$ and $\bfeta$, where $\bfeta$ takes the value $\eta$ with probability one; see for example the arguments around (3.7) in \cite{RockafellarUryasev.13}, which leverage the convexity and real-valuedness of $\cR_0$. For each $\omega\in \Omega$, the random variable appearing as an argument to $\cR_0$ has
\begin{align*}
&\psi_0^\nu\Big(\bfxi^\nu(\omega),F^\nu\big(\bfxi^\nu(\omega)\big),\bbH\big(G^\nu\big(\bfxi^\nu(\omega)\big)\big)\Big) - \psi_0\big(\bar\xi,x(\bar\xi),y(\bar\xi)\big)\\
& = \psi_0^\nu\big(\bfxi^\nu(\omega),x(\bar\xi),y(\bar\xi)\big) - \psi_0\big(\bfxi^\nu(\omega),x(\bar\xi),y(\bar\xi)\big) + \psi_0\big(\bfxi^\nu(\omega),x(\bar\xi),y(\bar\xi)\big) - \psi_0\big(\bar\xi,x(\bar\xi),y(\bar\xi)\big)\\
& \leq \Big|\psi_0^\nu\big(\bfxi^\nu(\omega),x(\bar\xi),y(\bar\xi)\big) - \psi_0\big(\bfxi^\nu(\omega),x(\bar\xi),y(\bar\xi)\big)\Big| + \Big|\psi_0\big(\bfxi^\nu(\omega),x(\bar\xi),y(\bar\xi)\big) - \psi_0\big(\bar\xi,x(\bar\xi),y(\bar\xi)\big)\Big|\\
& \leq \sigma^\nu + \kappa_0\, \diam \,\Xi.
\end{align*}
Since $\cR_0$ is monotone, this implies that
\[
  \mathfrak{m}^\nu \leq \sigma^\nu + \kappa_0\, \diam \,\Xi \, + \inf_{(x,y)\in C} \psi_0(\bar\xi,x,y) + \gamma.
\]
The claimed inequality follows for $\xi = \bar\xi$ because $\gamma$ is arbitrary.

Finally, consider the case with $\inf_{(x,y)\in C} \psi_0(\bar\xi,x,y) = -\infty$. Then, there exists $(x(\bar\xi),y(\bar\xi))\in C$ such that
\[
\psi_0\big(\bar\xi,x(\bar\xi),y(\bar\xi)\big) \leq -1/\gamma.
\]
Repeating the construction of $(F^\nu,G^\nu)$, which again is feasible for (TP)$^\nu$ without risk constraints \eqref{eqn:riskconstr}, we obtain that
\[
  \mathfrak{m}^\nu \leq \cR_0\Big(\psi_0^\nu\big(\bfxi^\nu,F^\nu(\bfxi^\nu),\bbH\big(G^\nu(\bfxi^\nu)\big)\big) - \psi_0\big(\bar\xi,x(\bar\xi),y(\bar\xi)\big)\Big)  - 1/\gamma.
\]
Thus, $\mathfrak{m}^\nu = -\infty$ because  $\gamma$ is arbitrary and the first term on right-hand side is bounded from above by $\sigma^\nu + \kappa_0\, \diam \,\Xi$ as argued above. The claimed inequality holds for $\xi = \bar \xi$ in this case as well.\eop\\

The theorem furnishes a uniform lower bound on the minimum value of the actual problem (AP)$(\xi)$ for $\xi\in \Xi$.  The lower bound can be compared to the objective function value $\psi_0(\xi,F^\nu(\xi), \bbH(G^\nu(\xi)))$ produced by any decision rule $(F^\nu,G^\nu)$. The comparison can determine for which $\xi\in \Xi$ the decision rule $(F^\nu,G^\nu)$ is satisfactory and for which the actual problem needs to be solved directly.

Theorem \ref{thm:lower} encourages construction of training problems involving the worst-case risk measure and a large number of training points because these choices drive up $\mathfrak{m}^\nu$ and produce a better lower bound. However, this must be traded off against the computational cost caused by many training points and also the (potential) suboptimality of the decision rules constructed from solving such training problems. The lower bound is also improved if $\diam\,\Xi$ is small, which implicitly means that the training points $\{\bfxi^\nu(\omega), \omega\in \Omega\}$ are tightly concentrated.

We next establish guarantees about the level of suboptimality for three decision rules originating from a near-minimizer $(F^\nu,G^\nu)$ of (TP)$^\nu$. While the first decision rule is simply $(F^\nu,G^\nu)$, the two others, called MDR and AMDR, are modifications ensuring feasibility in (AP)($\xi$) for all $\xi\in \Xi$.

MDR, an abbreviation for {\em minimum decision rule}, prescribes for every $\xi\in \Xi$ the decision  
\begin{align}
(x^\nu,y^\nu) & =  \Big(F^\nu\big(\bfxi^\nu(\omega^\nu)\big), \,\bbH\big(G^\nu\big(\bfxi^\nu(\omega^\nu)\big)\big)\Big),\label{eqn:MDRdef}\\
& \mbox{ where } \omega^\nu \in \nargmin_{\omega\in \Omega} \psi_0^\nu\Big(\bfxi^\nu(\omega),F^\nu\big(\bfxi^\nu(\omega)\big),\bbH\big(G^\nu\big(\bfxi^\nu(\omega)\big)\big)\Big).\nonumber
\end{align}
The choice $(x^\nu,y^\nu)$ is the decision prescribed by $(F^\nu,G^\nu)$ at the training point producing the lowest value of $\psi_0^\nu$. An advantage of MDR over $(F^\nu,G^\nu)$ is its feasibility in (AP)$(\xi)$ regardless of $\xi\in\Xi$. For other constant decision rules in specific settings, we refer to \cite{BertsimasGoyalLu.15,AwasthiGoyalLu.19}. 

AMDR, which stands for {\em adaptive minimum decision rule}, adjusts its prescription based on $\xi\in \Xi$ as specified by $(\bar F, \bar G)$, where
\begin{align}
\bar F(\xi) & = F^\nu\Big(\bfxi^\nu\big(\bar\omega(\xi)\big)\Big) ~\mbox{ and } ~\bar G(\xi) = G^\nu\Big(\bfxi^\nu\big(\bar\omega(\xi)\big)\Big),\label{eqn:AMDRdef}\\     
& \mbox{ with } \bar\omega(\xi) \in \nargmin_{\omega\in \Omega} \psi_0\Big(\xi,F^\nu\big(\bfxi^\nu(\omega)\big),\bbH\big(G^\nu\big(\bfxi^\nu(\omega)\big)\big)\Big).\nonumber
\end{align}
We find that $(\bar F(\xi), \bar G(\xi))$ is set to the decision prescribed by $(F^\nu,G^\nu)$ at one of the training points, with the selection depending on $\xi$. The selection of the training point, specified by $\bar\omega(\xi)$, requires access to the objective function $\psi_0$. While this is often unproblematic, it could represent a hurdle when the passing from $\psi_0$ to its approximation $\psi_0^\nu$ is indeed caused by computational challenges associated with $\psi_0$. Regardless, AMDR is easy to interpret: it prescribes at $\xi$ a decision with the lowest cost $\psi_0(\xi,\cdot,\cdot)$ among those identified by $(F^\nu,G^\nu)$ at the training points. Thus, it relates to $k$-adaptivity; see, e.g., \cite{BertsimasCaramanis.10,HanasusantoKuhnWiesemann.15}. AMDR always prescribes feasible decisions for (AP)$(\xi)$ at any $\xi\in \Xi$.

\begin{theorem}{\rm (decision rule guarantees).}\label{thm:error} For a fixed $\nu\in\nats$, suppose that $\delta,\delta^\nu \in (\epsilon,\infty]$, $\tau^\nu\in [0,\infty)$, $(F^\nu, G^\nu)$ is a $\tau^\nu$-minimizers of $\mbox{(TP)}^\nu$ without the risk constraints \eqref{eqn:riskconstr}, and assumptions (a)-(d) of Theorem \ref{thm:lower} hold. Then, the following three assertions are true:
\begin{enumerate}[(i)]
\item  The decision $(x^\nu,y^\nu)$ prescribed by MDR in \eqref{eqn:MDRdef} has the suboptimality guarantee 
\[
(x^\nu, y^\nu) \in (2\sigma^\nu + \tau^\nu + 2\kappa_0 \, \diam\, \Xi)\mbox{-}\nargmin_{(x,y)\in C} \psi_0(\xi,x,y)~~~~\forall \xi\in \Xi.
\]
\item  The AMDR, given by $(\bar F, \bar G)$ in \eqref{eqn:AMDRdef}, satisfies 
\[
\Big(\bar F(\xi), \bbH\big(\bar G(\xi)\big) \Big) \in (2\sigma^\nu + \tau^\nu + 2\kappa_0 \, \diam\, \Xi)\mbox{-}\nargmin_{(x,y)\in C} \psi_0(\xi,x,y)~~~~\forall \xi\in \Xi.
\]
\item The decision rule $(F^\nu, G^\nu)$ satisfies
\[
\Big(F^\nu(\xi), \bbH\big(G^\nu(\xi)\big) \Big) \in \Big(2\sigma^\nu + \tau^\nu + \big(\kappa_0+\kappa_0'\sqrt{1+\lambda^2}\big) \, \diam\, \Xi\Big)\mbox{-}\nargmin_{(x,y)\in C} \psi_0(\xi,x,y)
\]
for those $\xi\in \Xi$ producing $(F^\nu(\xi), \bbH(G^\nu(\xi))) \in C$ and $\bbH(G^\nu(\xi)) = \bbH(G^\nu(\bfxi^\nu(\omega)))$ for all $\omega\in \Omega$, where $\lambda$ and $\kappa_0'$ are given by
\begin{align*}
&\big\|F^\nu(\xi'') - F^\nu(\xi')\big\|_2 \leq \lambda \|\xi''-\xi'\|_2 ~~~\forall \xi'',\xi'\in \Xi\\
&\big|\psi_0(\xi'',x''  , y) - \psi_0(\xi',x'  , y)\big| \leq \kappa_0'\sqrt{\|\xi''-\xi'\|_2^2 + \|x''-x'\|_2^2}~~~\forall \xi'',\xi'\in \Xi, (x'',y)\in C, (x',y)\in C.
\end{align*}
\end{enumerate}
\end{theorem}
\state Proof. Let $\mathfrak{m}^\nu$ be the minimum value of $(\rm{TP})^\nu$ without the risk constraints \eqref{eqn:riskconstr}. By Theorem \ref{thm:lower}, one has for each $\xi\in \Xi$ that
\begin{equation}\label{eqn:errorproof}
\cR_0\Big(\psi_0^\nu\big(\bfxi^\nu,F^\nu(\bfxi^\nu),\bbH\big(G^\nu(\bfxi^\nu)\big)\big)\Big) + h^\nu(F^\nu,G^\nu) \leq \mathfrak{m}^\nu + \tau^\nu \leq \inf_{(x,y)\in C} \psi_0(\xi,x,y) + \sigma^\nu + \tau^\nu + \kappa_0 \, \diam\, \Xi.
\end{equation}
Consider (i) and let $\xi\in \Xi$ be arbitrary. Since $(F^\nu(\bfxi^\nu(\omega^\nu)), \bbH(G^\nu(\bfxi^\nu(\omega^\nu)))) \in C$, we also have $(x^\nu, y^\nu)\in C$. Thus, $\inf_{(x,y)\in C} \psi_0(\xi,x,y)\in \reals$; see also \eqref{eqn:errorproof}. We invoke \eqref{eqn:errorproof} to obtain
\begin{align}\label{eqn:prooferror0}
  & \psi_0(\xi,x^\nu, y^\nu) \leq  \psi_0(\xi,x^\nu, y^\nu)\\
  &~~~ + \inf_{(x,y)\in C} \psi_0(\xi,x,y) + \sigma^\nu + \tau^\nu + \kappa_0 \, \diam\, \Xi - \cR_0\Big(\psi_0^\nu\big(\bfxi^\nu,F^\nu(\bfxi^\nu),\bbH\big(G^\nu(\bfxi^\nu)\big)\big)\Big) - h^\nu(F^\nu,G^\nu)\nonumber\\
  & \leq  \inf_{(x,y)\in C} \psi_0(\xi,x,y) + \sigma^\nu + \tau^\nu + \kappa_0 \, \diam\, \Xi - \cR_0\Big(\psi_0^\nu\big(\bfxi^\nu,F^\nu(\bfxi^\nu),\bbH\big(G^\nu(\bfxi^\nu)\big)\big)-\psi_0(\xi,x^\nu, y^\nu)\Big),\nonumber
\end{align}
where we use the facts that $h^\nu(F^\nu,G^\nu)\geq 0$ and that $\cR_0(\bfeta' + \bfeta) = \cR_0(\bfeta') + \eta$ for any random variables $\bfeta'$ and $\bfeta$, where $\bfeta$ takes the value $\eta$ with probability one; see for example the arguments around (3.7) in \cite{RockafellarUryasev.13}, which leverage the convexity and real-valuedness of $\cR_0$. For each $\omega\in \Omega$, we obtain that
\begin{align*}
&\psi_0^\nu\Big(\bfxi^\nu(\omega),F^\nu\big(\bfxi^\nu(\omega)\big),\bbH\big(G^\nu\big(\bfxi^\nu(\omega)\big)\big)\Big)-\psi_0(\xi,x^\nu, y^\nu)\\
& = \psi_0^\nu\Big(\bfxi^\nu(\omega),F^\nu\big(\bfxi^\nu(\omega)\big),\bbH\big(G^\nu\big(\bfxi^\nu(\omega)\big)\big)\Big)-\psi_0^\nu\big(\bfxi^\nu(\omega^\nu),x^\nu,y^\nu\big) + \psi_0^\nu\big(\bfxi^\nu(\omega^\nu),x^\nu,y^\nu\big) - \psi_0(\xi,x^\nu, y^\nu)\\
& \geq \psi_0^\nu\big(\bfxi^\nu(\omega^\nu),x^\nu,y^\nu\big) - \psi_0\big(\bfxi^\nu(\omega^\nu),x^\nu,y^\nu\big) + \psi_0\big(\bfxi^\nu(\omega^\nu),x^\nu,y^\nu\big)- \psi_0(\xi,x^\nu, y^\nu)\\
& \geq -\sigma^\nu - \kappa_0\,\diam\,\Xi,
\end{align*}
where the second-to-last inequality follows by construction of $\omega^\nu$ and the last inequality from assumption (c) in Theorem \ref{thm:lower}. Since $\cR_0$ is monotone, this implies that
\[
\cR_0\Big(\psi_0^\nu\big(\bfxi^\nu,F^\nu(\bfxi^\nu),\bbH\big(G^\nu(\bfxi^\nu)\big)\big)-\psi_0\big(\xi,x^\nu, y^\nu)\Big) \geq -\sigma^\nu - \kappa_0\,\diam\,\Xi.
\]
Returning to \eqref{eqn:prooferror0}, we conclude that
\[
\psi_0(\xi,x^\nu,y^\nu) \leq \inf_{(x,y)\in C} \psi_0(\xi,x,y) + 2\sigma^\nu + \tau^\nu + 2\kappa_0 \, \diam\, \Xi.
\]

Consider (ii) and let $\xi\in \Xi$ be arbitrary. Since $(F^\nu(\bfxi^\nu(\bar\omega(\xi))), \bbH(G^\nu(\bfxi^\nu(\bar\omega(\xi))))) \in C$, we also have $(\bar F(\xi), \bbH(\bar G(\xi)))\in C$. Parallel to \eqref{eqn:prooferror0}, we invoke \eqref{eqn:errorproof} to establish
\begin{align}\label{eqn:prooferror01}
  & \psi_0\Big(\xi,\bar F(\xi), \bbH\big(\bar G(\xi)\big)\Big)\\
  & \leq  \inf_{(x,y)\in C} \psi_0(\xi,x,y) + \sigma^\nu + \tau^\nu + \kappa_0 \, \diam\, \Xi - \cR_0\Big(\psi_0^\nu\big(\bfxi^\nu,F^\nu(\bfxi^\nu),\bbH\big(G^\nu(\bfxi^\nu)\big)\big)-\psi_0\big(\xi,\bar F(\xi), \bbH\big(\bar G(\xi)\big)\big)\Big).\nonumber
\end{align}
For each $\omega\in \Omega$, we obtain that
\begin{align*}
&\psi_0^\nu\Big(\bfxi^\nu(\omega),F^\nu\big(\bfxi^\nu(\omega)\big),\bbH\big(G^\nu\big(\bfxi^\nu(\omega)\big)\big)\Big)-\psi_0\Big(\xi,\bar F(\xi), \bbH\big(\bar G(\xi)\big)\Big)\\
& = \psi_0^\nu\Big(\bfxi^\nu(\omega),F^\nu\big(\bfxi^\nu(\omega)\big),\bbH\big(G^\nu\big(\bfxi^\nu(\omega)\big)\big)\Big)-\psi_0\Big(\xi,F^\nu\big(\bfxi^\nu(\omega)\big),\bbH\big(G^\nu\big(\bfxi^\nu(\omega)\big)\big)\Big)\\
&~~~~~ + \psi_0\Big(\xi,F^\nu\big(\bfxi^\nu(\omega)\big),\bbH\big(G^\nu\big(\bfxi^\nu(\omega)\big)\big)\Big) - \psi_0\Big(\xi,F^\nu\big(\bfxi^\nu(\bar\omega(\xi))\big), \bbH\big(G^\nu\big(\bfxi^\nu(\bar\omega(\xi))\big)\big)\Big)\\
& \geq \psi_0^\nu\Big(\bfxi^\nu(\omega),F^\nu\big(\bfxi^\nu(\omega)\big),\bbH\big(G^\nu\big(\bfxi^\nu(\omega)\big)\big)\Big) - \psi_0\Big(\bfxi^\nu(\omega),F^\nu\big(\bfxi^\nu(\omega)\big),\bbH\big(G^\nu\big(\bfxi^\nu(\omega)\big)\big)\Big)\\
&~~~~~  + \psi_0\Big(\bfxi^\nu(\omega),F^\nu\big(\bfxi^\nu(\omega)\big),\bbH\big(G^\nu\big(\bfxi^\nu(\omega)\big)\big)\Big)  - \psi_0\Big(\xi,F^\nu\big(\bfxi^\nu(\omega)\big), \bbH\big(G^\nu\big(\bfxi^\nu(\omega)\big)\big)\Big)\\
& \geq -\sigma^\nu - \kappa_0\,\diam\,\Xi,
\end{align*}
where the second-to-last inequality follows by construction of $\bar\omega(\xi)$ and the last inequality from assumption (c) in Theorem \ref{thm:lower}. Since $\cR_0$ is monotone, this implies that
\[
\cR_0\Big(\psi_0^\nu\big(\bfxi^\nu,F^\nu(\bfxi^\nu),\bbH\big(G^\nu(\bfxi^\nu)\big)\big)-\psi_0\big(\xi,\bar F(\xi), \bbH\big(\bar G(\xi)\big)\big)\Big) \geq -\sigma^\nu - \kappa_0\,\diam\,\Xi.
\]
Returning to \eqref{eqn:prooferror01}, we conclude that
\[
\psi_0\Big(\xi,\bar F(\xi), \bbH\big(\bar G(\xi)\big)\Big) \leq \inf_{(x,y)\in C} \psi_0(\xi,x,y) + 2\sigma^\nu + \tau^\nu + 2\kappa_0 \, \diam\, \Xi.
\]
This implies the second assertion because $\inf_{(x,y)\in C} \psi_0(\xi,x,y)\in \reals$ as argued above.

Consider (iii) and let $\xi\in \Xi$ be such that $(F^\nu(\xi), \bbH(G^\nu(\xi))) \in C$ and $\bbH(G^\nu(\xi)) = \bbH(G^\nu(\bfxi^\nu(\omega)))$ for all $\omega\in \Omega$. Parallel to \eqref{eqn:prooferror0}, we invoke \eqref{eqn:errorproof} to establish
\begin{align}\label{eqn:prooferror02}
  & \psi_0\Big(\xi,F^\nu(\xi), \bbH\big(G^\nu(\xi)\big)\Big)\leq\\
  & \inf_{(x,y)\in C} \psi_0(\xi,x,y) + \sigma^\nu + \tau^\nu + \kappa_0 \, \diam\, \Xi - \cR_0\Big(\psi_0^\nu\big(\bfxi^\nu,F^\nu(\bfxi^\nu),\bbH\big(G^\nu(\bfxi^\nu)\big)\big)-\psi_0\Big(\xi,F^\nu(\xi), \bbH\big(G^\nu(\xi)\big)\Big)\Big).\nonumber
\end{align}
For each $\omega\in \Omega$, we obtain that
\begin{align*}
&\psi_0^\nu\Big(\bfxi^\nu(\omega),F^\nu\big(\bfxi^\nu(\omega)\big),\bbH\big(G^\nu\big(\bfxi^\nu(\omega)\big)\big)\Big)-\psi_0\Big(\xi,F^\nu(\xi), \bbH\big(G^\nu(\xi)\big)\Big)\\
& = \psi_0^\nu\Big(\bfxi^\nu(\omega),F^\nu\big(\bfxi^\nu(\omega)\big),\bbH\big(G^\nu\big(\bfxi^\nu(\omega)\big)\big)\Big)-\psi_0\Big(\bfxi^\nu(\omega),F^\nu\big(\bfxi^\nu(\omega)\big),\bbH\big(G^\nu\big(\bfxi^\nu(\omega)\big)\big)\Big)\\
&~~~~~ + \psi_0\Big(\bfxi^\nu(\omega),F^\nu\big(\bfxi^\nu(\omega)\big),\bbH\big(G^\nu\big(\bfxi^\nu(\omega)\big)\big)\Big) - \psi_0\Big(\xi,F^\nu(\xi), \bbH\big(G^\nu(\xi)\big)\Big)\\
& \geq -\sigma^\nu  - \kappa_0'\sqrt{\big\|\bfxi^\nu(\omega) - \xi\big\|_2^2 + \big\|F^\nu\big(\bfxi^\nu(\omega)\big) - F^\nu(\xi) \big\|_2^2}
\geq -\sigma^\nu - \kappa_0'\sqrt{1+\lambda^2}\,\diam\,\Xi,
\end{align*}
where the second-to-last inequality follows by assumption (c) in Theorem \ref{thm:lower}. Since $\cR_0$ is monotone, this implies that
\[
\cR_0\Big(\psi_0^\nu\big(\bfxi^\nu,F^\nu(\bfxi^\nu),\bbH\big(G^\nu(\bfxi^\nu)\big)\big)-\psi_0\big(\xi, F^\nu(\xi), \bbH\big(G^\nu(\xi)\big)\big)\Big) \geq -\sigma^\nu - \kappa_0'\sqrt{1+\lambda^2}\,\diam\,\Xi.
\]
Returning to \eqref{eqn:prooferror02}, we conclude that
\[
\psi_0\Big(\xi,F^\nu(\xi), \bbH\big(G^\nu(\xi)\big)\Big) \leq \inf_{(x,y)\in C} \psi_0(\xi,x,y) + 2\sigma^\nu + \tau^\nu + \big(\kappa_0+\kappa_0'\sqrt{1+\lambda^2}\big) \, \diam\, \Xi.
\]
This implies the third assertion because $\inf_{(x,y)\in C} \psi_0(\xi,x,y)\in \reals$ as argued above.\eop

Part (i) of the theorem examines MDR, which prescribes $(x^\nu,y^\nu)$ regardless of $\xi\in \Xi$.  While MDR is constant, it has a uniform suboptimality guarantee in (AP)$(\xi)$ across $\xi\in \Xi$. Part (ii) addresses AMDR, specified by $(\bar F, \bar G)$, which typically is {\em not} constant on $\Xi$. AMDR always prescribes feasible decisions for (AP)$(\xi)$, regardless of $\xi\in \Xi$, with the same suboptimality guarantee as for MDR. In part (iii), we consider $(F^\nu,G^\nu)$ directly. Since it may not prescribe feasible decisions for each $\xi\in \Xi$, we are unable to provide a uniform suboptimality guarantee for all $\xi\in \Xi$. However, if $(F^\nu,G^\nu)$ prescribes a feasible decision at $\xi$, which presumably can be checked easily after the solution of the training problem (TP)$^\nu$, then we obtain a suboptimality guarantee relative to (AP)$(\xi)$ provided that we also have $\bbH(G^\nu(\xi)) = \bbH(G^\nu(\bfxi^\nu(\omega)))$ for all $\omega\in \Omega$. The latter requirement is nontrivial but can be enforced a priori by selecting the mappings in $\cG$ to have a Lipschitz modulus $\kappa_\cG$ and $\epsilon$ in (TP)$^\nu$ to satisfy $\kappa_\cG \, \diam \, \Xi < \epsilon$. Specifically, if $\|G(\xi'') - G(\xi')\|_\infty \leq \kappa_\cG \|\xi''-\xi'\|_2$ for $\xi'',\xi'\in \Xi$ and $G\in \cG$, then the constraint \eqref{eqn:Dcon} ensures that each component function $g_i^\nu$ of $G^\nu$ satisfies either $g_i^\nu(\xi) > 0$ for all $\xi\in \Xi$ or $g_i^\nu(\xi) < 0$ for all $\xi\in \Xi$. If $\cG$ is further restricted to the affine mappings, then the requirement on $\epsilon$ can be relaxed to having $\kappa_\cG \, \diam \, \Xi < 2\epsilon$. These ``large'' values of $\epsilon$ must be introduced cautiously because they effectively restrict the attention to constant mappings in $\cG$, which could be meaningful in some applications but too costly in others. Nevertheless, the restriction is computationally beneficial; for example in (TP-super), $B$ and $b$ can be removed  and $y_i(\omega)$ can be replaced by $y_i$. If $F^\nu(\xi') = A^\nu\xi' + a^\nu$, then the parameter $\lambda$ in theorem becomes $\lambda = \|A^\nu\|$ using the matrix norm induced by $\|\cdot\|_2$.

For all three parts the assumptions are mild: expectation, superquantiles, worst-case, and many other risk measures can be used; the cost function can be nonconvex and nonsmooth; and the decision rules can be nearly arbitrary continuous mappings.

In summary among the three decision rules of Theorem \ref{thm:error}, MDR and AMDR have the same suboptimality guarantee while direct use of $(F^\nu,G^\nu)$ is theoretically inferior due to possible infeasibility and its worst suboptimality guarantee. We therefore recommend the use of MDR and AMDR in practice. The former has the advantage of being constant, which might be appealing to decision makers concerned about ``stability.'' One can expect the latter to prescribe better decisions due to its adaptivity, and this is also confirmed empirically in the next section.

\section{Decomposition Algorithm and Computational Results}\label{sec:num}

The broad framework of Sections \ref{sec:framework}-\ref{sec:error} applies in numerous settings. To illustrate some possible results and identify computational challenges, we consider a model from search theory motivated by \cite{RoysetSato.10,LejeuneRoysetMa.23}. 

A single searcher is looking for two moving targets in a discretized environment consisting of cells $c = 1, \dots, C$. In each time period $t \in \{1, \dots, T\}$, the searcher is in a cell and can look once in that cell for each target. The probability that a look finds a target is $g\in (0,1)$, given that the target is in the cell at that time, which is the glimpse detection probability. This probability is conveniently expressed using the detection rate $\alpha = - \ln(1-g)$. The targets' locations are random: In scenario $i \in \{1, \dots, I\}$, which occurs with probability $q_i \in [0,1]$, target $k$'s location at time $t$ is given by the binary vector $(\zeta_{1,t,i}^k, \dots, \zeta_{C,t,i}^k)$. Each binary vector has a single 1 specifying the location of the target at time $t$: $\zeta_{c,t,i}^k = 1$ if target $k$ is in cell $c$ at time $t$ and $\zeta_{c,t,i}^k = 0$ otherwise. Naturally, $\sum_{i=1}^I q_i = 1$. If the searcher is in cell $c$ at time $t<T$, then it can move to any cell in $N(c) \subset \{1, \dots, C\}$ for time period $t+1$. Let $y_{c,t} = 1$ if the searcher is in cell $c$ at time $t$ and $y_{c,t} = 0$ otherwise. The vector $y = (y_{1,1}, \dots, y_{C,T})$ is thus to be optimized. The scenario probabilities and detection rate depend on a varying parameter vector $\xi\in \reals^r$ and we write $q_i(\xi)$ and $\alpha(\xi)$ in place of $q_i$ and $\alpha$.

Under an independence assumption (see \cite{RoysetSato.10,LejeuneRoysetMa.23} for details), the problem of minimizing the nondetection probability for target 1 while keeping the nondetection probability for target 2 no higher than $\tau$ can be formulated as the binary convex program: 
\[
\mbox{(SP2)($\xi$)}~~~~~~  \nnmin_y ~\sum_{i = 1}^I  q_i(\xi) \exp\bigg( - \alpha(\xi) \sum_{c = 1}^C \sum_{t = 1}^T \zeta_{c,t,i}^1 \, y_{c,t}  \bigg)
\]
\begin{align}
  \mbox{subject to } \sum_{i = 1}^I q_i(\xi) \exp\bigg( - \alpha(\xi) &\sum_{c = 1}^C \sum_{t = 1}^T \zeta_{c,t,i}^2 \, y_{c,t}  \bigg) \leq \tau\label{eqn:secTargetCon}\\
  \sum_{c=1}^C y_{c,t} = 1~~~~&\forall t = 1, \dots, T\label{eqn:pathcon1}\\
                     \sum_{c'\in N(c)} y_{c',t-1} \geq y_{c,t} ~~~~&\forall c = 1, \dots, C, ~t=2, \dots, T\label{eqn:pathcon2}\\
                     y_{c,t}  \in \{0,1\}~~~~&\forall c=1, \dots, C,~ t=1, \dots, T.\nonumber
\end{align}
The program is a special case of the actual problem (AP)($\xi$) without continuous variables; we here concentrate on finding an affine decision rule for the more challenging case of binary variables. Thus, the decision rule is of the form $G(\xi) = B\xi + b$, where the $CT$-by-$r$ matrix $B$ and the $CT$-dimensional vector $b$ are to be determined by a training problem. The recommended decision for (SP2)($\xi$) becomes $y = \bbH(B\xi + b)$.

Throughout, we consider $T = 8$ time periods and $C = 81$ cells in a grid of 9-by-9 cells. Target 1 starts in cell 41 in the center and target 2 in cell 67; cell numbers begin in the upper-left corner and increase as we move right and down. (Figure \ref{fig:grid} in the Supplementary Material of Section \ref{sec:supp} shows cell numbers.) For time $t = 2, 3, \dots, T$, the targets are located according to $I = 100$ scenarios. (The scenarios are randomly sampled independently from a Markov chain in which a target remains in its current cell with probability 0.6 and moves to any available adjacent cell above, below, to the right, and to the left, with equal probability.) The set $N(c)$, constraining the searcher's movement (cf. \eqref{eqn:pathcon2}), contains $c$ and the cells above, below, to the right, and to the left of $c$, if available.

Even with the simplification of (SP2)($\xi$) having no continuous variables, a resulting training problem (TP)$^\nu$ tends to be large-scaled even for a moderate number of training points. This is further complicated by the nonlinearity in (SP2)($\xi$), which is best handled through a linearization approach and the addition of a large number of variables and constraints. (The Supplementary Material, Subsection \ref{subsec:implementationDetails}, provides details.) Risk measures such as those based on superquantiles may also prompt reformulations with auxiliary variables and constraints; see (TP-super) in Subsection \ref{subsec:implement} and Subsection \ref{subsec:implementationDetails}. A comprehensive study of how to best solve (TP)$^\nu$ in various settings is better postponed to a later paper. Here, we present a decomposition algorithm for training problems with a monotone risk measure in the objective function and worst-case risk measures in the constraints.

\subsection{Decomposition Algorithm}\label{DCM}

For the purpose of constructing a relaxation, we can remove the nonnegative term $h^\nu(F,G)$ from (TP)$^\nu$ and also replace $F(\bfxi^\nu(\omega))$ and $\bbH(G(\bfxi^\nu(\omega)))$ by $\omega$-dependent decisions $x(\omega) \in \reals^n$ and $y(\omega) \in \{0,1\}^m$, respectively. It then becomes apparent that solving (AP)$(\bfxi^\nu(\omega))$ for each $\omega\in\Omega$ leads to a lower bound on the minimum value of (TP)$^\nu$ under the assumption that $\cR_0$ is monotone and $\cR_1, \dots, \cR_q$ are worst-case risk measures. Case-dependent approaches can furnish an upper bound. We describe the details for instances with linearly independent training points  $\{\bfxi^\nu(\omega) \in \reals^r, \omega\in \Omega\}$, $\cG$ consisting of affine mappings of the form $B\xi + b$, and a regularization term $h^\nu(F,G) = h(F) + \theta \sum_{i=1}^m \|B_i\|$ for some $\theta\in [0,\infty)$, $h:\cF\to [0,\infty)$, and $G(\xi) = (g_1(\xi), \dots, g_m(\xi))$ with $g_i(\xi) = \langle B_i, \xi\rangle + b_i$. The vector norm $\|\cdot\|$ is arbitrary.\\

\vspace{0.2cm}

\state Decomposition Algorithm for Training Problem (TP)$^\nu$. 

\begin{description}

\item[Step 1.] For each $\omega\in \Omega$, compute a feasible solution $(\bfx^\nu(\omega),   \bfy^\nu(\omega))$ of (AP)$(\bfxi^\nu(\omega))$ and a lower bound ${\bfL}^\nu(\omega)$ on its minimum value. Set $L^\nu = \cR_0(\bfL^\nu)$. 

\item[Step 2.] For each $i=1, \dots, m$, set $  B^\nu_{i}$ and $  b^\nu_{i}$ as follows:

If $  \bfy^\nu_{i}(\omega) = 1$ for all $\omega \in \Omega$, set $  B^\nu_{i} = 0$ and $  b^\nu_{i} = \epsilon$.  

If $  \bfy^\nu_{i}(\omega) = 0$ for all $\omega \in \Omega$, set $  B^\nu_{i} = 0$ and $  b^\nu_{i} = -\epsilon$.  

Otherwise, solve the convex separation problem 
\begin{align*}
	\nnmin_{B_{i}, b_{i}} \,\|B_{i}\|  ~\mbox{ subject to } ~ \epsilon \leq \big\langle  B_{i}, \bfxi^\nu(\omega) \big\rangle + b_{i}  &\leq \delta^\nu  &\forall \omega \mbox{ with }   \bfy^\nu_{i}(\omega) = 1&\\  
	-\epsilon \geq   \big\langle  B_{i}, \bfxi^\nu(\omega) \big\rangle + b_{i}  &\geq -\delta^\nu &\forall \omega \mbox{ with }   \bfy^\nu_{i}(\omega) = 0&
\end{align*}
and let $(  B^\nu_{i},  b^\nu_{i})$ be a resulting minimizer.   

\item[Step 3.] Obtain a feasible solution $F^\nu$ of the problem: 
\begin{align}
\nnmin_{F\in \cF}  ~\cR_0\Big(\psi_0^\nu\big(\bfxi^\nu,F(\bfxi^\nu),  \bfy^\nu\big)\Big) + h(F) &+ \theta \sum_{i=1}^m \|  B^\nu_i\| &&\nonumber\\
 \mbox{subject to } ~~\psi_k^\nu\Big(\bfxi^\nu(\omega),F\big(\bfxi^\nu(\omega)\big),  \bfy^\nu\big(\bfxi^\nu(\omega)\big)\Big) & \leq 0, && k=1, \dots, q, ~\omega\in \Omega\nonumber\\
                           \Big(F\big(\bfxi^\nu(\omega)\big),   \bfy^\nu\big(\bfxi^\nu(\omega)\big)\Big) & \in C && \forall \omega\in \Omega.\nonumber
\end{align}
Let $U^\nu = \cR_0(\psi_0^\nu\big(\bfxi^\nu,F^\nu(\bfxi^\nu),  \bfy^\nu)) + h(F^\nu) + \theta \sum_{i=1}^m \|  B^\nu_i\|$. Return the decision rule given by $F^\nu$, $B^{\nu}$, and $b^{\nu}$, which is a feasible solution of (TP)$^\nu$ with relative optimality gap $(U^{\nu} - L^\nu)/L^\nu$. 

\end{description}

We observe that $L^\nu$ is a valid lower bound on the minimum value of (TP)$^\nu$ because $\cR_0$ is monotone. Step 2 produces $(B^\nu, b^\nu)$ that minimizes the regularization term $\theta \sum_{i=1}^m \|B_i\|$ while satisfying $B\bfxi^\nu(\omega) + b \in D_\epsilon(\delta^\nu)$ and $\bbH(B\bfxi^\nu(\omega) + b) = \bfy^\nu(\omega)$ for all $\omega\in \Omega$. This is computable using convex optimization when $\delta^\nu$ is sufficiently large and $\{\bfxi^\nu(\omega) \in \reals^r, \omega\in \Omega\}$ are linearly independent, which can be anticipated because $r$ is often larger than the number of training points. In our implementation, we use $\|\cdot\| = \|\cdot\|_1$ and solve Step 2 using linear programming after a standard linearization; see for example \cite[Section 2.H]{primer}. Any feasible solution $F^\nu$ from Step 3 produces together with $(B^\nu,b^\nu)$ a feasible decision rule in (TP)$^\nu$ because \eqref{eqn:Dcon} is satisfied by construction in Step 2. If no feasible solution exists in Step 3 or the optimality gap is too large, then we recommend applying a heuristic to (TP)$^\nu$ to obtain the upper bound $U^\nu$, which again can be compared with the lower bound $L^\nu$. In our numerical experiments, we solve a restriction of (TP)$^\nu$ with fixed $B_{i} = 0$ and $b_{i} = -\epsilon$ if $\bfy^\nu_{i}(\omega) = 0$ for all $\omega \in \Omega$ and warm-start the remaining $(B_i,b_i)$ at the value $(B^\nu_i,b^\nu_i)$.  

 Since Step 1 of the Decomposition Algorithms involves solving (AP)$(\bfxi^\nu(\omega))$ for each $\omega\in \Omega$, its computational cost is significantly higher than solving a single actual problem (AP)$(\xi)$. Some of this increase, however, might be avoided by allowing an optimality tolerance in Step 1 and by warm-starting the computations for (AP)$(\bfxi^\nu(\omega))$ by a solution for (AP)$(\bfxi^\nu(\omega'))$, with $\omega\neq \omega'$.

In the following, we test the Decomposition Algorithm in the context of the search problem $\mbox{(SP2)($\xi$)}$ and leverage Gurobi (version 10.0.1) as solver for the resulting linear programs and linear binary programs (using default options). We run the calculations on an Optiflex7090x60 machine with IntelCore i5-1055 3.10GHz processor and 12 logical processors. The subproblems in Steps 1 and 2 are solved to optimality. Step 3 is terminated after 300 seconds or when the optimality gap of the restricted problem fall below 0.01\%, whichever occurs first. Throughout, we use the parameter values $\tau=0.45$, $\theta = 0.001$, $\delta^\nu = 1$, and $\epsilon = 0.001$. The training points are assigned the same probability, i.e., $P(\omega) = 1/|\Omega|$.

\subsection{Examples of Convergence and Infeasibility}

We illustrate Theorem \ref{thm:convergence} and challenges related to infeasibility by considering a parameter vector $\xi = (\xi_0, \xi_1, \dots, \xi_{100}) \in \reals^r$, with $r = 101$, which influences the detection rate and scenario probabilities in (SP2)($\xi$) as follow: 
\[
\alpha(\xi) = \bar\alpha + \xi_0; ~~~~q_i(\xi)  = \frac{\max\{0, 1/100 + \xi_i\}}{\sum_{i'=1}^{100} \max\{0, 1/100 + \xi_{i'}\}}, ~~~i=1, \dots, 100,
\]
where $\bar\alpha = 2.74887$ is the nominal detection rate (also used in \cite{LejeuneRoysetMa.23} and corresponds to a glimpse detection probability  $g = 0.936$). Let $|\Omega| = 100$, i.e., there are 100 training points. We generate eight training data sets $\{\bfxi^\nu(\omega), \omega\in \Omega\}$, $\nu = 1, \dots, 8$. The first one, $\{\bfxi^1(\omega), \omega\in \Omega\}$, consists of points ``far'' from the zero vector and the eighth set, $\{\bfxi^8(\omega), \omega\in \Omega\}$, consists of points ``near'' the zero vector to illustrate the convergence $\bfxi^\nu \to \bar \bfxi$ as $\nu$ increases in Theorem  \ref{thm:convergence}, with $\bar \bfxi$ being the constant random vector with the zero vector as value. For each $\nu=1, \dots, 8$ and $\omega \in \Omega$, we generate $\bfxi^\nu(\omega) = (\bfxi^\nu_0(\omega), \bfxi^\nu_1(\omega), \dots, \bfxi^\nu_{100}(\omega))$ 
by setting $\bfxi^\nu_0(\omega)$ equal to a number sampled from a uniform distribution on $[-0.18 + 0.02\nu, ~0.18 - 0.02\nu]$ and by setting $\bfxi^\nu_i(\omega)$ equal to a number sampled from a uniform distribution on $[-0.009 + 0.001\nu, ~0.009 - 0.001\nu]$. (All the sampling assumes statistical independence.)

We consider two versions of the risk measure $\cR_0$: expectation and worst-case, with $\cR_1, \dots, \cR_q$ being worst-case risk measures in both cases. These versions produce two instances of (TP)$^\nu$ called $\mbox{(EW-SP2)}^\nu$ and $\mbox{(WW-SP2)}^\nu$. The Supplementary Material in Subsection \ref{subsec:implementationDetails} provides detailed formulations of these training problems.

For each $\nu = 1, \dots, 8$, we solve the resulting (EW-SP2)$^\nu$ using the Decomposition Algorithm and obtain $(B^\nu,b^\nu)$. Table \ref{DEC2} summarizes the results. The first row labeled $\nu = 1$ relies on the training points $\{\bfxi^1(\omega), \omega\in \Omega\}$ and reports the prescribed decision $\bbH(B^1 \cdot 0 + b^1)$ for (SP2)(0). Columns 2-9 give the specific cells to search for each of the eight time periods. It is immediately clear that the decision rule fails to produce a feasible path. In several time periods, the decision rule prescribes the searcher to be in multiple cells at the same time. This violates constraint \eqref{eqn:pathcon1}. Columns 10 and 11 give the nondetection probability for target 1 and target 2, respectively, for the prescribed decision (despite its infeasibility). The last column reports a 22\% relative optimality gap ($100\%(U^{\nu} - L^\nu)/L^\nu$) from the Decomposition Algorithm.

\begin{table}[ht]
		\centering
		\setlength\extrarowheight{1pt}
		\footnotesize{
			\begin{tabular}{c|c|c|c|c|c|c|c|c|c|c|c}
				\hline
				\multirow{2}{*}{$\nu$} 
				&\multicolumn{8}{c|}{searcher cells across time}
				&prob.
				&prob.
				&gap\\ 
				\cline{2-9} 
				& 1 & 2 &3&4&5&6&7&8 & tar. 1 & tar. 2 & (\%)\\
				\hline
				$1$ &67&40,58&49&40&31,40,41,49&32,40,58&40,49,66&40,57,75 &0.287  & 0.029& 22\\ 			
				$2$&67&40,58&49&50,58&41,49&40,41,50,67,68&40,41,66,67&40,75 &0.278  &0.022& 26\\		
				$3$&67&58&49&50,58&31&40,76&41,67&40 &0.422 &0.038& 24\\ 
				$4$&67&58&40&-& 31&67&40,66&40,75,76&0.504 &0.033& 21\\
				$5$&67&58&49&68&31,41&40,41,67&41,66&40,75,76 &0.424 &0.026& 28\\ 
				$6$&67&40,58&49&50&31,41,49&40,41,50&41,77&40  &0.323 &0.048& 29\\ 
				$7$&67&58&49&40&31&40&41,77&40,76 &0.398 &0.050& 25\\ 
				$8$&{\bf 67}&{\bf 58}&{\bf 49}&{\bf 40}&{\bf 31}&{\bf 40}&{\bf 41}&{\bf 40} &0.490 &0.057& 0.0\\ 
				$\infty$&{\bf 67}&{\bf 58}&{\bf 49}&{\bf 40}&{\bf 31}&{\bf 40}&{\bf 41}&{\bf 40}& 0.490 &0.057 & 0\\
				\hline
			\end{tabular}
			\caption{Prescribed cells using (EW-SP2)$^\nu$ and $|\Omega| = 100$ training points. Boldface indicates a sequence of cells that satisfies the constraints \eqref{eqn:pathcon1} and \eqref{eqn:pathcon2}. The dash specifies no search in that time period.}\label{DEC2}
		}
\end{table}

For comparison, the last row of the table gives a minimizer of (SP2)(0) and the corresponding probabilities of nondetection. Theorem \ref{thm:convergence} asserts that under certain assumptions (which are satisfied for the present instance) the prescribed decision $\bbH(b^\nu)$ will eventually match a minimizer of (SP2)(0) as $\nu\to \infty$. Table \ref{DEC2} confirms this fact with the second-to-last row (using training data $\{\bfxi^8(\omega), \omega\in \Omega\}$) indeed matching the last row.

\begin{table}[ht]
		\centering
		\setlength\extrarowheight{1pt}
		\footnotesize{
			\begin{tabular}{c|c|c|c|c|c|c|c|c|c|c|c}
				\hline
				\multirow{2}{*}{$\nu$} 
				&\multicolumn{8}{c|}{searcher cells across time}
				&prob.
				&prob.
				&gap\\ 
				\cline{2-9} 
				& 1 & 2 &3&4&5&6&7&8 & tar. 1 & tar. 2 & (\%)\\
				\hline
				$1$&67&58&49&40&31,41,49&32,40,58,67&40,66&40,75 &0.362 &0.030& 10\\ 			
				$2$&67&58&49&50&41&41&41,49&58 &0.483 &0.051& 2.7\\
				$3$&67&58&49&50,58&31&40,50,76&41,67&40,58,66 &0.408 & 0.057&12\\						
				$4$&{\bf 67}&{\bf 58}&{\bf 49}&{\bf 40}&{\bf 41}&{\bf 41}&{\bf 41}&{\bf 40} &0.491 &0.057& 3.4\\ 
				$5$&{\bf 67}&{\bf 58}&{\bf 49}&{\bf 50}&{\bf 41}&{\bf 41}&{\bf 41}&{\bf 40} &0.491 &0.057& 1.6\\ 
				$6$&{\bf 67}&{\bf 58}&{\bf 49}&{\bf 40}&{\bf 31}&{\bf 40}&{\bf 41}&{\bf 40} &0.490 &0.057& 3.7\\ 
				$7$&{\bf 67}&{\bf 58}&{\bf 49}&{\bf 40}&{\bf 31}&{\bf 40}&{\bf 41}&{\bf 40} &0.490 &0.057& 3.7\\ 
				$8$&{\bf 67}&{\bf 58}&{\bf 49}&{\bf 50}&{\bf 41}&{\bf 41}&{\bf 41}&{\bf 40} &0.491 &0.057 &0.1\\ 
				$\infty$&{\bf 67}&{\bf 58}&{\bf 49}&{\bf 40}&{\bf 31}&{\bf 40}&{\bf 41}&{\bf 40}& 0.490 & 0.057 & 0\\
				\hline
			\end{tabular}
			\caption{Prescribed cells using (WW-SP2)$^\nu$ and $|\Omega| = 100$ training points. Boldface indicates a sequence of cells that satisfies the constraints \eqref{eqn:pathcon1} and \eqref{eqn:pathcon2}.}
			\label{DEC4}
		}
\end{table}

We replicate the effort but now by solving the training problem (WW-SP2)$^\nu$, i.e., relying on the worst-case risk measure for $\cR_0$. Table \ref{DEC4} summarizes the results obtained by the Decomposition Algorithm. In this case, the convergence to a minimizer of (SP2)(0) seems to take place quicker with $\nu = 6$ and $7$ already replicating the last row in the table. The rows labeled $\nu = 4, 5, 8$ only marginally miss the minimizer; the prescribed decisions give a nondetection probability of 0.491, which is only 0.001 above the minimum value of (SP2)(0).  

The improved results might be caused by the fact that the Decomposition Algorithm performs better on (WW-SP2)$^\nu$ than on (EW-SP2)$^\nu$: optimality gaps are reduced from 22\% to 5\%, on average; see the last columns in Tables \ref{DEC2} and \ref{DEC4}. In either case, the computing times are moderate, around 15 minutes, as detailed in Table \ref{runtimes} in the Supplementary Material.

Tables \ref{DC1} and \ref{DEC3} replicate the effort in Tables \ref{DEC2} and \ref{DEC4} but now using only 10 training points, i.e., $|\Omega| = 10$. Again (WW-SP2)$^\nu$ tends to produce the better decision rules compared to (EW-SP2)$^\nu$. Less training data makes it easier to achieve decisions that satisfy the constraints \eqref{eqn:pathcon1} and \eqref{eqn:pathcon2}, but the decisions might still be suboptimal; see for example the rows labeled $\nu = 3,5,6$ in Table \ref{DC1}. In these cases, the decision rules prescribe searching first in cell 41, which is the initial location of target 1. While this brings down dramatically the nondetection probability for that target, it makes the decision fall short of meeting the 0.45-probability constraint \eqref{eqn:secTargetCon} for target 2. Thus, the decisions in the rows labeled $\nu = 3,5,6$ are infeasible for (SP2)(0). The minimizer of (SP2)(0) in the last row in the tables entails searching first in cell 67, where target 2 starts. This guarantees feasibility with respect to the 0.45-probability constraint \eqref{eqn:secTargetCon}. A training problem undoubtedly will have difficulty with picking up on these subtleties, but as we see from the tables it eventually manages as the training points approach the zero vector. In fact, this is guaranteed by Theorem \ref{thm:convergence}.

\begin{table}[ht]
	\centering
	\setlength\extrarowheight{1pt}
	\footnotesize{
			\begin{tabular}{c|c|c|c|c|c|c|c|c|c|c|c}
				\hline
				\multirow{2}{*}{$\nu$} 
				&\multicolumn{8}{c|}{searcher cells across time}
				&prob.
				&prob.
				&gap\\ 
				\cline{2-9} 
				& 1 & 2 &3&4&5&6&7&8 & tar. 1 & tar. 2 & (\%)\\
				\hline
			$1$&41&40&49&58&67&50,67,68&66,67&66,75 &0.032&0.343& 3.6\\
			$2$&67&58&49&40,50&41,49&32,40,41&41&40 &0.414&0.056& 2.1\\ 
			$3$&{\bf 41}&{\bf 50}&{\bf 59}&{\bf 58}&{\bf 67}&{\bf 68}&{\bf 67}&{\bf 66}&0.045&0.456& 3.3\\ 
			$4$&41&50&59&40,58&31,67&40,68&67&66   &0.033&0.456& 3.5\\ 
			$5$&{\bf 41}&{\bf 50}&{\bf 59}&{\bf 58}&{\bf 67}&{\bf 68}&{\bf 67}&{\bf 66}&0.045&0.456& 2.9\\ 			
			$6$&{\bf 41}&{\bf 50}&{\bf 59}&{\bf 58}&{\bf 67}&{\bf 68}&{\bf 67}&{\bf 66}&0.045&0.456& 3.1\\ 
			$7$&{\bf 67}&{\bf 58}&{\bf 49}&{\bf 40}&{\bf 49}&{\bf 50}&{\bf 41}&{\bf 40}&0.491&0.056& 0.2\\ 
			$8$&{\bf 67}&{\bf 58}&{\bf 49}&{\bf 40}&{\bf 31}&{\bf 40}&{\bf 41}&{\bf 40}&0.490& 0.057& 0.1\\ 
			$\infty$&{\bf 67}&{\bf 58}&{\bf 49}&{\bf 40}&{\bf 31}&{\bf 40}&{\bf 41}&{\bf 40}&0.490&0.057& 0\\
			\hline
		\end{tabular}
		\caption{Prescribed cells using (EW-SP2)$^\nu$ and $|\Omega| = 10$ training points. Boldface indicates a sequence of cells that satisfies the constraints \eqref{eqn:pathcon1} and \eqref{eqn:pathcon2}.} \label{DC1}
	}
\end{table}

\begin{table}[ht]
	\centering
	\setlength\extrarowheight{1pt}
	\footnotesize{
			\begin{tabular}{c|c|c|c|c|c|c|c|c|c|c|c}
				\hline
				\multirow{2}{*}{$\nu$} 
				&\multicolumn{8}{c|}{searcher cells across time}
				&prob.
				&prob.
				&gap\\ 
				\cline{2-9} 
				& 1 & 2 &3&4&5&6&7&8 & tar. 1 & tar. 2 & (\%)\\
				\hline
			$1$&41&50&59&58,68&67&50,67,68&66,77&75,76                                  &0.044 &0.241 &1.4\\	
			$2$&{\bf 67}&{\bf 58}&{\bf 49}&{\bf 50}&{\bf 41}&{\bf 41}&{\bf 41}&{\bf 40} &0.491 &0.057 &0.0\\ 	
			$3$&{\bf 67}&{\bf 58}&{\bf 49}&{\bf 40}&{\bf 31}&{\bf 40}&{\bf 41}&{\bf 40} &0.490 &0.057 &0.0\\ 
			$4$&{\bf 67}&{\bf 58}&{\bf 49}&{\bf 40}&{\bf 31}&{\bf 40}&{\bf 41}&{\bf 40} &0.490 &0.057 &0.1\\ 
			$5$&{\bf 67}&{\bf 58}&{\bf 49}&{\bf 40}&{\bf 49}&{\bf 50}&{\bf 41}&{\bf 40} &0.491 &0.056 &0.0\\
			$6$&{\bf 67}&{\bf 58}&{\bf 49}&{\bf 40}&{\bf 41}&{\bf 41}&{\bf 41}&{\bf 40} &0.491 &0.057 &0.1\\ 
			$7$&{\bf 67}&{\bf 58}&{\bf 49}&{\bf 40}&{\bf 31}&{\bf 40}&{\bf 41}&{\bf 40} &0.490 &0.057 &0.1\\ 
			$8$&{\bf 67}&{\bf 58}&{\bf 49}&{\bf 40}&{\bf 49}&{\bf 50}&{\bf 41}&{\bf 40} &0.495 &0.056 &0.0\\ 
			$\infty$&{\bf 67}&{\bf 58}&{\bf 49}&{\bf 40}&{\bf 31}&{\bf 40}&{\bf 41}&{\bf 40} &0.490 &0.057 & 0\\
			\hline
		\end{tabular}
	\caption{Prescribed cells using (WW-SP2)$^\nu$ and $|\Omega| = 10$ training points. Boldface indicates a sequence of cells that satisfies the constraints \eqref{eqn:pathcon1} and \eqref{eqn:pathcon2}.}
	\label{DEC3}
}	
\end{table}

Generally, fewer training points imply less computing time for the Decomposition Algorithm; on average less than one minute as seen from Table \ref{runtimes} in the Supplementary Material.

An extreme choice would be to adopt a single training point $\xi^\nu$ (and then $|\Omega| = 1$). While Theorem \ref{thm:convergence} applies even in this case, the resulting decision rule becomes trivial: prescribe a minimizer of (SP2)($\xi^\nu$) for any other instance (SP2)($\xi$) including when $\xi= 0$. (Note that (EW-SP2)$^\nu$ and (WW-SP2)$^\nu$ are equivalent in this case and it suffices to solve (SP2)($\xi^\nu$) to obtain the decision rule.) Simulations not reported on in detail here show convergence as well as the trivial fact that all prescribed decisions satisfy the constraints \eqref{eqn:pathcon1} and \eqref{eqn:pathcon2}. However, a single training point cannot be expected to result in a ``robust'' decision rule that holds up for varying $\xi$. 

%
%

In summary, the empirical results confirm Theorem \ref{thm:convergence} and show that the training problems indeed prescribe optimal decisions for the actual problem when the training points are concentrated. This takes place for both expectation and worst-case risk measures. Still, as can be expected for a difficult combinatorial optimization such as (SP2)($\xi$), direct application of a solution of a training problem may produce infeasible decisions in the actual problem when the training points are dispersed. The next subsection addresses this concern.

\subsection{Examples of Feasible Decision Rules}\label{subsec:feasRules}

We now turn to the two decision rules in Theorem \ref{thm:error}, MDR and AMDR, which are guaranteed to produce feasible decisions regardless of the number of training points. We also quantify the frequency by which a decision rule {\em directly} obtained from a training problem is feasible. Since Theorem \ref{thm:error} omits risk constraints, we consider in this subsection instances with a single target, i.e., (SP2)($\xi$) without the constraint \eqref{eqn:secTargetCon}, which we refer to as (SP1)($\xi$). We focus on the instance of (TP)$^\nu$ obtained by adopting 
the $\beta$-superquantile, $\beta\in (0,1)$, as the risk measure $\cR_0$; a detailed formulation of the resulting training problem referred to as $\mbox{($\beta$-SP1)}^\nu$ appears in the Supplementary Material (see Subsection \ref{subsec:implementationDetails}). Subsection \ref{subsub:addnum} furnishes additional numerical results for expectation and worst-case risk measures. The resulting training problems are referred to as  (E-SP1)$^\nu$ and (W-SP1)$^\nu$, respectively.

We also switch to a parameter vector $\xi = (\xi_1, \dots, \xi_{100}) \in \reals^{100}$, which influences the scenario probabilities according to the expressions $q_i(\xi) = 1/100 + \xi_i$, $i=1, \dots, 100$. The detection rate is $\alpha(\xi) = \bar\alpha = 0.510826$ for all $\xi$, which corresponds to a glimpse detection probability  $g = 0.4$. These adjustments recalibrate the instances so that the minimum value of (SP1)$(\xi)$ is approximately 0.4 and allow us to easily compare with the error bounds in Theorem \ref{thm:error}. 

We consider 100 training points (each being a 100-dimensional vector): $\bfxi^\nu(\omega) = (\bfxi^\nu_1(\omega), \dots, \bfxi^\nu_{100}(\omega))$ $\in$ $\reals^{100}$, $\omega \in \Omega$ so $|\Omega| = 100$. Each of these training points are generated by random sampling 100 times from a uniform distribution on $[0,1]$, normalizing the values so they sum to one and thus become a probability vector, and subtracting the vector $(1/100, \dots, 1/100)$. The resulting vector becomes a training point if its length is no larger than 0.05. We refer to these training points as the {\em uniform training data}.


\begin{table}[ht]
	\centering
	\setlength\extrarowheight{1pt}
	\footnotesize{
		\begin{tabular}{c|c|c|c|c|c|c}
			\hline
			\multirow{2}{*}{rule} 
			& training & test	&number 
			&\multicolumn{3}{c}{suboptimality}\\ 
			\cline{5-7} 
			& data & data & feasible          & min  &avg   & max\\
			\hline
			$B^\nu,b^\nu$ & unif & unif & 98  &0.000 &0.007 &0.022\\	
			MDR           & unif & unif & 100 &0.001 &0.011 &0.024\\	
			AMDR          & unif & unif & 100 &0.000 &0.002 &0.013\\				
			\hline
			$B^\nu,b^\nu$ & unif & beta & 99  &0.000 &0.014 &0.039\\	
			MDR           & unif & beta & 100 &0.000 &0.017 &0.045\\	
			AMDR          & unif & beta & 100 &0.000 &0.008 &0.033\\				
			\hline
			$B^\nu,b^\nu$ & beta & beta & 29  &0.001 &0.012 &0.039\\	
			MDR           & beta & beta & 100 &0.004 &0.023 &0.049\\	
			AMDR          & beta & beta & 100 &0.000 &0.001 &0.014\\
        \hline
		\end{tabular}
		\caption{Performance of decision rules obtain from ($\beta$-SP1)$^\nu$. Suboptimality statistics are computed over feasible decisions.}
		\label{tab:testS}
	}	
\end{table}

We produce a decision rule $(B^\nu,b^\nu)$ by solving ($\beta$-SP1)$^\nu$ using the Decomposition Algorithm, which takes 351 seconds. The resulting relative optimality gap is 3.62\%, which corresponds to an absolute gap of 0.014. (In this test and the following ones, we increase the 300-second limit in Step 3 of the Decomposition Algorithm to 1000 seconds.) The quality of the decision rule is assessed across 100 test points $\xi^k\in \reals^{100}$, $k=1, \dots, 100$, generated exactly as the uniform training data. We refer to these test points as the {\em uniform test data}. The first row of Table \ref{tab:testS} summarizes the results. The column labeled ``number feasible'' gives the number of test points $\xi^k$ at which $\bbH(B^\nu \xi^k + b^\nu)$ is feasible in (SP1)($\xi^k$), i.e., satisfies constraints \eqref{eqn:pathcon1} and \eqref{eqn:pathcon2}. We also solve (SP1)($\xi^k$) and obtain its minimum value $\mathfrak{m}(\xi^k)$. The level of suboptimality for $\bbH(B^\nu \xi^k + b^\nu)$ relative to (SP1)($\xi^k$) is quantified by computing
\[
\sum_{i = 1}^I  q_i(\xi^k) \exp\bigg( - \bar\alpha \sum_{c = 1}^C \sum_{t = 1}^T \zeta_{c,t,i}^1 \, \bbH\big( \langle B^\nu_{c,t}, \xi^k\rangle + b^\nu_{c,t} \big) \bigg) - \mathfrak{m}(\xi^k). 
\]
The last three columns in Table \ref{tab:testS} provide the minimum, average, and maximum values of these levels of suboptimality across the uniform test data. We see that the decision rule $(B^\nu,b^\nu)$ prescribes a feasible decision 98 out of 100 times. Across the feasible decisions, the level of suboptimality is 0.007, on average, which is rather minor when the minimum values $\mathfrak{m}(\xi^k)$, $k=1, \dots, 100$ are approximately 0.4.

The next two rows in Table \ref{tab:testS} provide parallel results for the two decision rules MDR and AMDR in Theorem \ref{thm:error}. Regardless of $\xi$, MDR prescribes the decision  
\[
y^\nu = \bbH\big(B^\nu \bfxi^\nu(\omega^\nu) + b^\nu\big), ~~\omega^\nu \in \nargmin_{\omega\in \Omega}  \sum_{i = 1}^I  q_i\big(\bfxi^\nu(\omega)\big) \exp\bigg( - \bar\alpha \sum_{c = 1}^C \sum_{t = 1}^T \zeta_{c,t,i}^1 \, \bbH\big( \langle B^\nu_{c,t}, \bfxi^\nu(\omega)\rangle + b^\nu_{c,t} \big) \bigg). 
\]  
AMDR prescribes the $\xi$-adaptive decision 
\[
\bar y(\xi) = \bbH\big(B^\nu \bfxi^\nu\big(\bar\omega(\xi)\big) + b^\nu\big), ~~ \bar\omega(\xi) \in \nargmin_{\omega\in \Omega}  \sum_{i = 1}^I  q_i(\xi) \exp\bigg( - \bar\alpha \sum_{c = 1}^C \sum_{t = 1}^T \zeta_{c,t,i}^1 \, \bbH\big( \langle B^\nu_{c,t}, \bfxi^\nu(\omega)\rangle + b^\nu_{c,t} \big) \bigg). 
\]
Thus, MDR and AMDR are both easily derived from $(B^\nu,b^\nu)$. By construction, MDR and AMDR always prescribe feasible decisions as seen in column 4. 
The levels of suboptimality for MDR are slightly worst than those from direct use of  $(B^\nu,b^\nu)$, but the latter decision rule is infeasible twice. AMDR is better than the two others; its worst test point exhibits a suboptimality of only 0.013. 

We also consider test points that are generated by random sampling 100 times from a beta distribution (parameters 0.1 and 0.1), normalizing the values so they sum to one and thus become a probability vector, and subtracting the vector $(1/100, \dots, 1/100)$. The resulting vector becomes a test point if its length is no larger than 0.1. We refer to these test points as the {\em beta test data}; they are more spread out than the uniform training data. Rows 4-6 in Table \ref{tab:testS} report results for the decision rules from rows 1-3 (obtained using the uniform training data) on the beta test data. While the levels of suboptimality are larger, the relative performance of the decision rules remains the same with AMDR providing high-quality recommendations.

To further validate the decision rules, we generate training data using the beta distribution as described above; we call it the {\em beta training data}. Again, we solve ($\beta$-SP1)$^\nu$ using the Decomposition Algorithm, which takes 125 seconds. The resulting relative optimality gap is 4.41\%, which corresponds to an absolute gap of 0.017. This produces a new decision rule $(B^\nu,b^\nu)$ and thus also new instances of MDR and AMDR. The last three rows of Table \ref{tab:testS} report the performance of these rules tested on the beta test data. (Additional results with the uniform test data are similar; see Table \ref{tab:testSextra}.) Since the beta training and test data are more dispersed than their uniform counterparts, direct use of $(B^\nu,b^\nu)$ results in infeasible decisions for around two-thirds of the test points. In contrast, AMDR remains highly effective prescribing near-minimizers consistently. We note that adaptation to the varying $\xi$ is important, with AMDR outperforming the constant decision rule MDR with a wide margin. 

We end by considering a more challenging target. Above, target 1 starts in cell 67 and fans out from there according to 100 scenarios. We now let it start in any of the nine cells in the middle row of the 9-by-9 grid, with equal probability, and then it proceeds according to 100 scenarios generated as before. We call this the {\em dispersed target} case. Table \ref{tab:testSdifferent} reports on the results using the same format as Table \ref{tab:testS}. The Decomposition Algorithm computes $(B^\nu, b^\nu)$ in 546 seconds in this case on the uniform training data; the relative optimality gap is 4.58\% and the absolute gap is 0.026. For the beta training data, the parallel numbers are 132 second, 4.59\%, and 0.026.  Again, direct use of $(B^\nu,b^\nu)$ produces infeasible decisions on more than half the test points. AMDR prescribe feasible decisions with objective function values within 0.04 of the actual minimum values across all test points examined in Table \ref{tab:testSdifferent}. The average value is an order of magnitude less and thus corresponds to a relative suboptimality gap of about 1\%. In comparison, MDR prescribes significantly worst decisions. We conclude that AMDR consistently is superior and exhibits significant robustness to training with {\em incorrectly distributed} data. (Additional results with the uniform test data appear in Table \ref{tab:testSdifferentextra}.)


\begin{table}[ht]
	\centering
	\setlength\extrarowheight{1pt}
	\footnotesize{
		\begin{tabular}{c|c|c|c|c|c|c}
			\hline
			\multirow{2}{*}{rule} 
			& training & test	&number 
			&\multicolumn{3}{c}{suboptimality}\\ 
			\cline{5-7} 
			& data & data & feasible          & min  &avg   & max\\
			\hline
			$B^\nu,b^\nu$ & unif & unif & 60  &0.000 &0.005 &0.043\\	
			MDR           & unif & unif & 100 &0.000 &0.003 &0.027\\	
			AMDR          & unif & unif & 100 &0.000 &0.002 &0.020\\				
			\hline
			$B^\nu,b^\nu$ & unif & beta & 39  &0.000 &0.027 &0.077\\	
			MDR           & unif & beta & 100 &0.000 &0.018 &0.075\\	
			AMDR          & unif & beta & 100 &0.000 &0.012 &0.037\\				
			\hline
			$B^\nu,b^\nu$ & beta & beta & 34  &0.003 &0.054 &0.135\\	
			MDR           & beta & beta & 100 &0.025 &0.069 &0.200\\	
			AMDR          & beta & beta & 100 &0.000 &0.006 &0.025\\
        \hline
		\end{tabular}
		\caption{Performance of decision rules obtain from ($\beta$-SP1)$^\nu$ against dispersed target. Suboptimality statistics are computed over feasible decisions.}
		\label{tab:testSdifferent}
	}	
\end{table}

Figure \ref{fig:plot} compares the performance of AMDR across all test instances and all training problems. The blue, orange, and grey bars in Figure \ref{fig:plot}(left) report average suboptimality for AMDR based on training problems (E-SP1)$^\nu$, ($\beta$-SP1)$^\nu$, and (W-SP1)$^\nu$, respectively. (Detailed results for (E-SP1)$^\nu$ and (W-SP1)$^\nu$,  and additional ones for ($\beta$-SP1)$^\nu$, appear in Subsection \ref{subsub:addnum}.) The corresponding bars in Figure \ref{fig:plot}(right) give the more conservative 95\% quantiles for the suboptimality values. Almost uniformly, ($\beta$-SP1)$^\nu$ is the more effective training problem for AMDR, with (W-SP1)$^\nu$ often only slightly worse. The expectation-based training problem  (E-SP1)$^\nu$ exhibits larger levels of suboptimality in most instances. This highlights the importance of experimenting with different risk measures beyond expected values. These findings are consistent with result from machine learning and stochastic optimization (see, e.g., \cite{LevyCarmonDuchiSidford.20,LaguelPillutlaMalickHarchaoui.21b,LaguelPillutlaMalickHarchaoui.21,GotohKimLim.21}), where some degree of risk-averseness (distributional robustness) during training tends to improve accuracy in out-of-sample tests.

\drawing{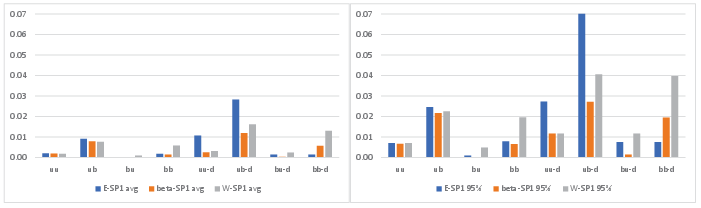}{6.8in} {Summary of suboptimality for AMDR across training problems and instances: uu = uniform training data and uniform test data, ub = uniform training data and beta test data, etc.; ``-d'' indicates the dispersed target case. Left portion gives average suboptimality and right portion gives 95\% quantile of suboptimality levels.}{fig:plot}

The level of suboptimality observed in Tables \ref{tab:testS} and \ref{tab:testSdifferent} and Figure \ref{fig:plot} are smaller than those predicted by Theorem \ref{thm:error}. In the notation of that theorem, we have $\sigma^\nu = 0$ because function values are not approximated. The tolerance $\tau^\nu$ in the solution of the training problems (by the Decomposition Algorithm) tends to be around 0.02 as reported above. For all $\xi,\xi'\in \Xi$ and $y\in \{0, 1\}^{CT}$ satisfying \eqref{eqn:pathcon1} and \eqref{eqn:pathcon2}, one has
\[
\Bigg|\sum_{i = 1}^I  q_i(\xi) \exp\bigg( - \bar\alpha \sum_{c = 1}^C \sum_{t = 1}^T \zeta_{c,t,i}^1 \, y_{c,t}\bigg) - \sum_{i = 1}^I  q_i(\xi') \exp\bigg( - \bar\alpha\sum_{c = 1}^C \sum_{t = 1}^T \zeta_{c,t,i}^1 \, y_{c,t}\bigg)\Bigg| \leq \sqrt{I}\|\xi- \xi'\|_2.
\] 
Thus, the Lipschitz modulus $\kappa_0 = \sqrt{I} = 10$ because the exponential terms are between 0 and 1. For the uniform test data and beta test data, $\diam\, \Xi$ equals $0.05$ and $0.10$, respectively. Thus, the level of suboptimality according to Theorem \ref{thm:error} is bounded by $2\sigma^\nu + \tau^\nu + 2\kappa_0 \,\diam \,\Xi = 0.02 + 2\cdot 10 \cdot 0.05 = 1.02$ for the uniform test data and twice that number for the beta test data. 
Tables \ref{tab:testS} and \ref{tab:testSdifferent} report significantly smaller numbers. One can expect, however, that sampling more than 100 test points will narrow this gap between theory and practical performance. The above calculations also utilize a ``conservative'' $\kappa_0$. If $y$-vectors under consideration had been narrowed, then $\kappa_0$ could have been taken as a lower number.

\section{Conclusions}

The decision rules developed in this paper help managers recognize relationships between data and decisions, and prescribe courses of actions that are guaranteed near-optimal as well as transparent and interpretable. Transparency and easy interpretability are prerequisites for any decision support tool that needs to be conveyed, understood, and justified broadly. An understanding of the relationship between data and decisions is important for managers as they identify the critical factors, prioritize their focus, guide the resource allocation, and develop mitigation strategies for handling data uncertainty. 

The relationship between data and decision can be especially convoluted for combinatorial optimization problems. We consider simple decisions rules (for example, affine ones) that ensure transparency and interpretability, but also allow nearly any continuous decision rule within a unified framework of analysis. We also permit a wide range of risk measures that promote robust and conservative decision rules. These measures include risk averse ones that exhibit strong empirical performance as compared to the use of expected values. Our preferred Adaptive Minimum Decision Rule (AMDR) prescribes decisions with 1\% relative suboptimality gap, on average, even when facing distributional shifts in out-of-sample testing.\\

\state Acknowledgement. This work of the first author is supported in part by ONR (N000142412277). The second author acknowledges partial support  from NSF (ECCS-2114100; RISE-2220626) and ONR (N00014-22-1-2649).

\bibliographystyle{plain}
\bibliography{refs}

\newpage

\section{Supplementary Material}\label{sec:supp}

\subsection{Proofs}

\state Proof of Proposition \ref{prop:recovery}.  Consider (a). Let $(F,G) \in \cF\times\cG$ satisfy \eqref{eqn:trainingproblemCon}-\eqref{eqn:trainingproblemYcon}. Since $\cR_k$ is the worst-case risk measure, one has for each $\omega\in \Omega$
\[
  \psi_k\Big(\bfxi(\omega),F\big(\bfxi(\omega)\big),\bbH\big(G\big(\bfxi(\omega)\big)\big)\Big) \leq 0, ~~k=1, \dots, q.
\]
Since $(\bar F(\bfxi(\omega)),\bbH(\bar G(\bfxi(\omega))))$ is a minimizer of (AP)$(\bfxi(\omega))$ and $(F(\bfxi(\omega)),\bbH(G(\bfxi(\omega))))$ is feasible in the same problem, one has $(\bar F(\bfxi(\omega)),\bbH(\bar G(\bfxi(\omega))))\in C$ and
\begin{align*}
  \psi_0\Big(\bfxi(\omega),\bar F\big(\bfxi(\omega)\big),\bbH\big(\bar G\big(\bfxi(\omega)\big)\big)\Big) &\leq \psi_0\Big(\bfxi(\omega),F\big(\bfxi(\omega)\big),\bbH\big(G\big(\bfxi(\omega)\big)\big)\Big)\\
  \psi_k\Big(\bfxi(\omega),\bar F\big(\bfxi(\omega)\big),\bbH\big(\bar G\big(\bfxi(\omega)\big)\big)\Big) & \leq 0, ~~k=1, \dots, q,
\end{align*}
with these conditions holding for each $\omega\in \Omega$. By monotonicity of $\cR_0, \cR_1, \dots, \cR_q$, this implies that
\begin{align*}
  \cR_0\Big(\psi_0\big(\bfxi,\bar F(\bfxi),\bbH\big(\bar G(\bfxi)\big)\big)\Big) &\leq \cR_0\Big(\psi_0\big(\bfxi,F(\bfxi),\bbH\big(G(\bfxi)\big)\big)\Big)\\
  \cR_k\Big(\psi_k\big(\bfxi,\bar F(\bfxi),\bbH\big(\bar G(\bfxi)\big)\big)\Big) & \leq 0, ~~k=1, \dots, q.
\end{align*}
Thus, $(\bar F,\bar G)$ is feasible in (TP) and, in fact, a minimizer because $(F,G)$ is chosen arbitrarily.

Next, consider (b). In view of the fact that $\cR_1, \dots, \cR_q$ are worst-case risk measures, the decision $(F^\star(\bfxi(\omega)), \bbH(G^\star(\bfxi(\omega))))$ is feasible for (AP)$(\bfxi(\omega))$ regardless of $\omega\in \Omega$. Since $(\bar F(\bfxi(\omega)), \bbH(\bar G(\bfxi(\omega))))$ is a minimizer of (AP)$(\bfxi(\omega))$, this implies that
\begin{equation}\label{eqn:proofrecovery}
\psi_0\Big(\bfxi(\omega),\bar F\big(\bfxi(\omega)\big),\bbH\big(\bar G\big(\bfxi(\omega)\big)\big)\Big) \leq \psi_0\Big(\bfxi(\omega),F^\star\big(\bfxi(\omega)\big),\bbH\big(G^\star\big(\bfxi(\omega)\big)\big)\Big) ~~~\forall \omega\in \Omega.
\end{equation}
We observe that $(\bar F, \bar G)$ is feasible in (TP) as argued above. This fact together with the optimality of $(F^\star,G^\star)$ in (TP), with $\cR_0$ being the expectation, imply that
\[
   \sum_{\omega\in \Omega} P(\omega) \psi_0\Big(\bfxi(\omega),F^\star\big(\bfxi(\omega)\big),\bbH\big(G^\star\big(\bfxi(\omega)\big)\big)\Big) \leq \sum_{\omega\in \Omega} P(\omega) \psi_0\Big(\bfxi(\omega),\bar F\big(\bfxi(\omega)\big),\bbH\big(\bar G\big(\bfxi(\omega)\big)\big)\Big).
   \]
Since $P(\omega)>0$ for all $\omega\in \Omega$, the inequalities in \eqref{eqn:proofrecovery} must hold with equality. We conclude that $(F^\star(\bfxi(\omega)), \bbH(G^\star(\bfxi(\omega))))$ is actually a minimizer of (AP)$(\bfxi(\omega))$ for each $\omega\in \Omega$.\eop\\

\state Proof of Theorem \ref{thm:convergence}. The argument leverages epi-convergence. We define $\phi^\nu,\phi:\cF\times\cG \to [-\infty,\infty]$ by setting
\begin{align}\label{eqn:phinudef}
\phi^\nu(F,G) & = \cR_0\Big(\psi_0^\nu\big(\bfxi^\nu,F(\bfxi^\nu),\bbH\big(G(\bfxi^\nu)\big)\big)\Big) + h^\nu(F,G) + \sum_{k=1}^q \iota_{(-\infty,0]} \Big(\cR_k\Big(\psi_k^\nu\big(\bfxi^\nu,F(\bfxi^\nu),\bbH\big(G(\bfxi^\nu)\big)\big)\Big)\Big)\nonumber\\
& ~~~~+ \sum_{\omega \in \Omega} \iota_C\Big(F\big(\bfxi^\nu(\omega)\big), \bbH\big(G\big(\bfxi^\nu(\omega)\big)\big)\Big) + \sum_{\omega \in \Omega} \iota_{D_\epsilon(\delta^\nu)}\Big( G\big(\bfxi^\nu(\omega)\big) \Big)\\
\phi(F,G) & = \cR_0\Big(\psi_0\big(\bar\bfxi,F(\bar\bfxi),\bbH\big(G(\bar\bfxi)\big)\big)\Big) + h(F,G) + \sum_{k=1}^q \iota_{(-\infty,0]} \Big(\cR_k\Big(\psi_k\big(\bar\bfxi,F(\bar\bfxi),\bbH\big(G(\bar\bfxi)\big)\big)\Big)\Big)\nonumber\\
& ~~~~+ \sum_{\omega \in \Omega} \iota_C\Big(F\big(\bar\bfxi(\omega)\big), \bbH\big(G\big(\bar\bfxi(\omega)\big)\big)\Big) + \sum_{\omega \in \Omega} \iota_{D_\epsilon(\delta)}\Big( G\big(\bar\bfxi(\omega)\big) \Big),\nonumber
\end{align}
where for any set $D$, $\iota_D(x) = 0$ when $x\in D$ and $\iota_D(x) = \infty$ otherwise. We recall that $\phi^\nu$ epi-converges to $\phi$ provided that (see, e.g., \cite{Royset.18})
\begin{align}
  &\forall (F^\nu,G^\nu) \to (F,G), ~\nliminf \phi^\nu(F^\nu,G^\nu) \geq \phi(F,G)\label{eqn:liminfcond}\\
  &\forall (F,G)\in \cF\times\cG, ~~\exists (F^\nu,G^\nu) \to (F,G) \, \mbox{ such that }  \, \nlimsup \phi^\nu(F^\nu,G^\nu) \leq \phi(F,G).\label{eqn:limsupcond}
\end{align}
We first establish \eqref{eqn:liminfcond}.  Suppose that $(F^\nu,G^\nu) \to (F,G)$. Thus, $\max_{\xi\in \Xi} \|F^\nu(\xi)-F(\xi)\|_2\to 0$. As a consequence of $\bfxi^\nu\to \bar\bfxi$, one has $\bfxi^\nu(\omega)\in \Xi\to \bar\bfxi(\omega)$ for each $\omega\in \Omega$, with the right-hand side being equal to $\bar\xi\in \Xi$. This implies that $F^\nu(\bfxi^\nu(\omega)) \to F(\bar\xi)$ for each $\omega\in \Omega$. A similar argument establishes that $G^\nu(\bfxi^\nu(\omega)) \to G(\bar\xi)$ for all $\omega\in \Omega$.

If $(F^\nu(\bfxi^\nu(\omega)), \bbH(G^\nu(\bfxi^\nu(\omega))))\not\in C$ for some $\omega\in \Omega$, then $\phi^\nu(F^\nu,G^\nu) = \infty$. Similarly, if $G^\nu(\bfxi^\nu(\omega))\not\in D_\epsilon(\delta^\nu)$ for some $\omega\in \Omega$, then $\phi^\nu(F^\nu,G^\nu) = \infty$. We can therefore assume without loss of generality that $(F^\nu(\bfxi^\nu(\omega)), \bbH(G^\nu(\bfxi^\nu(\omega))))\in C$ and $G^\nu(\bfxi^\nu(\omega))\in D_\epsilon(\delta^\nu)$ for every $\omega\in \Omega$ and $\nu\in\nats$ in verification of \eqref{eqn:liminfcond}. This in turn implies that $G(\bar\bfxi(\omega)) \in D_\epsilon(\delta)$ for every $\omega\in \Omega$ because $\delta^\nu\to \delta$. Hence, there is $\bar\nu$ such that $\bbH(G^\nu(\bfxi^\nu(\omega))) = \bbH(G(\bar\bfxi(\omega)))$ for all $\nu\geq \bar\nu$ and $\omega\in \Omega$. Since $F^\nu(\bfxi^\nu(\omega))\to F(\bar\bfxi(\omega))$ and $C$ is closed, we also have
$(F(\bar\bfxi(\omega)), \bbH(G(\bar\bfxi(\omega))))\in C$ for all $\omega\in \Omega$. Moreover, by assumption (b),
\[
\psi_k^\nu\Big(\bfxi^\nu(\omega),F^\nu\big(\bfxi^\nu(\omega)\big),\bbH\big(G\big(\bfxi^\nu(\omega)\big)\big)\Big)\to  \psi_k\Big(\bar\bfxi(\omega),F\big(\bar\bfxi(\omega)\big),\bbH\big(G\big(\bar\bfxi(\omega)\big)\big)\Big), ~~k=0, 1, \dots, q,
\]
and this holds for each $\omega\in \Omega$. Thus, the random variables $\psi_k^\nu(\bfxi^\nu,F^\nu(\bfxi^\nu),\bbH(G(\bfxi^\nu)))$ converge pointwise to $\psi_k(\bar\bfxi,F(\bar\bfxi),\bbH(G(\bar\bfxi)))$, $k=0, 1, \dots, q$. The continuity of $\cR_0, \cR_1, \dots, \cR_q$ implies that
\[
\cR_0\Big(\psi_0^\nu\big(\bfxi^\nu,F^\nu(\bfxi^\nu),\bbH\big(G(\bfxi^\nu)\big)\big)\Big)\to  \cR_0\Big(\psi_0\big(\bar\bfxi,F(\bar\bfxi),\bbH\big(G(\bar\bfxi)\big)\big)\Big)
\]
and, for $k=1, \dots, q$,
\[
\nliminf \iota_{(-\infty,0]}\Big(\cR_k\Big(\psi_k^\nu\big(\bfxi^\nu,F^\nu(\bfxi^\nu),\bbH\big(G(\bfxi^\nu)\big)\big)\Big)\Big)\geq \iota_{(-\infty,0]}\Big(\cR_k\Big(\psi_k\big(\bar\bfxi,F(\bar\bfxi),\bbH\big(G(\bar\bfxi)\big)\big)\Big)\Big).
\]
Since none of the terms in the sum defining $\phi(F,G)$ equals $-\infty$ and, by assumption (d), $\nliminf h^\nu(F^\nu,G^\nu)$ $\geq$ $h(F,G)$, we obtain that
\begin{align*}
\nliminf \phi^\nu(F^\nu,G^\nu) & \geq  \nliminf \cR_0\Big(\psi_0^\nu\big(\bfxi^\nu,F^\nu(\bfxi^\nu),\bbH\big(G^\nu(\bfxi^\nu)\big)\big)\Big) + \nliminf h^\nu(F^\nu,G^\nu)\\
& ~~~+ \sum_{k=1}^q \nliminf \iota_{(-\infty,0]} \Big(\cR_k\Big(\psi_k^\nu\big(\bfxi^\nu,F^\nu(\bfxi^\nu),\bbH\big(G^\nu(\bfxi^\nu)\big)\big)\Big)\Big)\\
& \geq \cR_0\Big(\psi_0\big(\bar\bfxi,F(\bar\bfxi),\bbH\big(G(\bar\bfxi)\big)\big)\Big) + h(F,G) + \sum_{k=1}^q\iota_{(-\infty,0]}\Big(\cR_k\Big(\psi_k\big(\bar\bfxi,F(\bar\bfxi),\bbH\big(G(\bar\bfxi)\big)\big)\Big)\Big)\\
& = \phi(F,G)
\end{align*}
and thus \eqref{eqn:liminfcond} holds.

We next establish \eqref{eqn:limsupcond}. Let $(F,G)\in \cF \times \cG$. Without loss of generality we assume that $G(\bar\bfxi(\omega)) \in D_\epsilon(\delta)$ and $(F(\bar\bfxi(\omega)), \bbH(G(\bar\bfxi(\omega)))) \in C$ for all $\omega\in \Omega$ and $\cR_k(\psi_k(\bar\bfxi,F(\bar\bfxi),\bbH(G(\bar\bfxi))))\leq 0$ for $k= 1, \dots, q$ because otherwise $\phi(F,G) = \infty$ and \eqref{eqn:limsupcond} holds trivially. For each $\nu\in\nats$, we construct $G^\nu:\Xi\to \reals^m$ by setting
\[
G^\nu(\xi) = G(\xi) + b^\nu, ~~\xi\in \Xi,
\]
where $b^\nu\in \reals^m$ is to be determined. By assumption (c), $G^\nu\in \cG$. We define $b^\nu = (b_1^\nu, \dots, b_m^\nu)$ componentwise and let
\[
  b_i^\nu = \begin{cases}
    \max_{\omega\in \Omega}\max\big\{0, -\delta^\nu - g_i(\bfxi^\nu(\omega))\big\}     & \mbox{ if } g_i(\bar \xi) = -\delta\\
    \min_{\omega\in \Omega}\min\big\{0, -\epsilon - g_i(\bfxi^\nu(\omega))\big\} & \mbox{ if } g_i(\bar \xi) = -\epsilon\\
    \max_{\omega\in \Omega}\max\big\{0, \epsilon - g_i(\bfxi^\nu(\omega))\big\} & \mbox{ if } g_i(\bar \xi) = \epsilon\\
    \min_{\omega\in \Omega}\min\big\{0, \delta^\nu - g_i(\bfxi^\nu(\omega))\big\} & \mbox{ if } g_i(\bar \xi) = \delta\\
    0 & \mbox{ otherwise,}
       \end{cases}
\]
where $g_i$ is the $i$th component function of $G$. Since $\delta^\nu\to \delta$ and $g_i(\bfxi^\nu(\omega)) \to g_i(\bar\xi)$ for every $\omega\in \Omega$, one has $b_i^\nu\to 0$. Thus, $\max_{\xi\in \Xi}\|G^\nu(\xi) - G(\xi)\|_2\to 0$. Moreover, for any $\omega \in \Omega$,
\[
  g_i^\nu\big(\bfxi^\nu(\omega)\big) = g_i\big(\bfxi^\nu(\omega)\big) + b_i^\nu \in \begin{cases}
    [-\delta^\nu,\infty)     & \mbox{ if } g_i(\bar \xi) = -\delta\\
    (-\infty,-\epsilon] & \mbox{ if } g_i(\bar \xi) = -\epsilon\\
    [\epsilon, \infty) & \mbox{ if } g_i(\bar \xi) = \epsilon\\
    (-\infty,\delta^\nu] & \mbox{ if } g_i(\bar \xi) = \delta,
       \end{cases}
\]
where $g_i^\nu$ is the $i$th component function of $G^\nu$. This means that $G^\nu(\bfxi^\nu(\omega))\in D_\epsilon(\delta^\nu)$ for sufficiently large $\nu$ regardless of $\omega\in \Omega$. In fact, there exists $\bar\nu$ such that $\bbH(G^\nu(\bfxi^\nu(\omega))) = \bbH(G(\bar\bfxi(\omega)))$ for all $\nu\geq \bar\nu$ and $\omega\in \Omega$.

Next, we turn to the construction of $F^\nu$. Set $\bar y = \bbH(G(\bar\xi))$. By assumption (e), there exist sequences $\{x^\mu \in \reals^n, \mu\in \nats\}$ and $\{\gamma^\mu \in (0,\infty), \mu\in \nats\}$ with the properties that $x^\mu\to F(\bar \xi)$, $\gamma^\mu \to 0$, and
\begin{align*}
\psi_k\big(\bar\xi,x^\mu,\bar y\big) &< 0, ~k=1, \dots, q, ~~~\mu\in\nats\\
(z^\mu, \bar y) &\in C ~~\forall z^\mu \in \ball(x^\mu,\gamma^\mu), ~~~\mu\in\nats.
\end{align*}
We define $\{\bar F^\mu, \mu\in\nats\}$ by setting
\[
\bar F^\mu(\xi) = F(\xi) - F(\bar\xi) + x^\mu ~~\mbox{ for } \xi\in \Xi.
\]
By assumption (c), $\bar F^\mu\in \cF$. By construction, $\max_{\xi\in \Xi} \|\bar F^\mu(\xi) - F(\xi)\|_\infty \to 0$ as $\mu\to \infty$ and, for each $\mu\in\nats$, one has
\[
\psi_k\big(\bar\xi,\bar F^\mu(\bar\xi),\bar y\big) < 0, ~k=1, \dots, q ~~~~\mbox{ and }~~~~ (z^\mu, \bar y) \in C ~~\forall z^\mu \in \ball\big(\bar F^\mu(\bar\xi),\gamma^\mu\big).
\]
Let $\bar\nu_0 = \bar\nu$. For each fixed $\mu\in\nats$, there exists $\bar\nu_\mu > \bar\nu_{\mu-1}$ such that for all $\nu\geq \bar\nu_\mu$ one has
\begin{align*}
  \psi_k^\nu\big(\bfxi^\nu(\omega),\bar F^\mu(\bfxi^\nu(\omega)),\bar y\big) &\leq 0~~~~\forall \omega \in \Omega, ~~k=1, \dots, q\\
  \Big(\bar F^\mu\big(\bfxi^\nu(\omega)\big), \bar y\Big) & \in C~~~\forall \omega \in \Omega.
\end{align*}
We recall that $\bar F^\mu(\bfxi^\nu(\omega))\to \bar F^\mu(\bar \xi)$ as $\nu\to \infty$. Thus, the inequality holds for sufficiently large $\nu$ because of assumption (b).  The inclusion in $C$ must hold for sufficiently large $\nu$ because $\|\bar F^\mu(\bfxi^\nu(\omega))-\bar F^\mu(\bar \xi)\|_2\leq \gamma^\mu$ for sufficiently large $\nu$. We then construct $\{F^\nu\}_{\nu = \bar\nu_1}^\infty$ by setting $F^\nu = \bar F^\mu$ for $\bar \nu_{\mu+1} > \nu \geq \bar\nu_\mu$, $\mu\in \nats$. By construction, $(F^\nu(\bfxi^\nu(\omega)), \bar y) \in C$ for all $\omega \in \Omega$ and $\nu \geq \bar\nu_1$ and $\max_{\xi\in \Xi} \|F^\nu(\xi) - F(\xi)\|_\infty \to 0$. Moreover, because $\cR_1, \dots, \cR_q$ are monotone by assumption (a), we find that for $\nu\geq \bar\nu_1$,
\[
\cR_k\Big(\psi_k^\nu\big(\bfxi^\nu,F^\nu(\bfxi^\nu),\bar y\big)\Big) \leq 0, ~~k=1, \dots, q.
\]
We have constructed $(F^\nu,G^\nu)\to (F,G)$ with the property that
\begin{equation}\label{eqn:limsuptemp}
\nlimsup \phi^\nu(F^\nu,G^\nu) \leq \nlimsup \cR_0\Big(\psi_0^\nu\big(\bfxi^\nu,F^\nu(\bfxi^\nu),\bbH\big(G^\nu(\bfxi^\nu)\big)\big)\Big) + \nlimsup h^\nu(F^\nu,G^\nu).
\end{equation}
Since $\psi_0^\nu(\bfxi^\nu(\omega),F^\nu(\bfxi^\nu(\omega)),\bar y)\to \psi_0(\bar\bfxi(\omega), F(\bar\bfxi(\omega)),\bar y)$ for all $\omega\in \Omega$ by assumption (b), we find that the random variables
\[
\psi_0^\nu\Big(\bfxi^\nu,F^\nu(\bfxi^\nu),\bbH\big(G^\nu(\bfxi^\nu)\big)\Big)\to \psi_0\Big(\bar\bfxi,F(\bar\bfxi),\bbH\big(G(\bar\bfxi)\big)\Big).
\]
The continuity of $\cR_0$ (cf. assumption (a)) ensures that
\[
\cR_0\Big(\psi_0^\nu\big(\bfxi^\nu,F^\nu(\bfxi^\nu),\bbH\big(G^\nu(\bfxi^\nu)\big)\big)\Big) \to \cR_0\Big(\psi_0\big(\bar\bfxi,F(\bar\bfxi),\bbH\big(G(\bar\bfxi)\big)\big)\Big).
\]
These facts together with assumption (d) establish that \eqref{eqn:limsuptemp} implies \eqref{eqn:limsupcond}. We have confirmed that $\phi^\nu$ epi-converges to $\phi$.

Next, let $N\subset\nats$ be a subsequence and $\{(\hat F^\nu, \hat G^\nu), \nu\in\nats\}$ be a sequence satisfying assumption (f) with $(\hat F^\nu, \hat G^\nu)\Nto (F^\star,G^\star)$. It is well known  that
\[
\phi^\nu(\hat F^\nu,\hat G^\nu) \Nto \inf_{(F,G)\in \cF\times\cG} \phi(F,G) ~~\mbox{ and } ~~ (F^\star,G^\star) \in \nargmin_{(F,G)\in \cF\times\cG} \phi(F,G)
\]
whenever $\inf_{(F,G)\in \cF\times\cG} \phi(F,G)<\infty$; see, e.g., \cite[Proposition 2.1]{Royset.18} or \cite[Theorem 5.5]{primer} for a proof that extends to the present setting nearly verbatim. Directly from the definition of risk measures, we find that for any $(F,G)\in \cF\times\cG$,
\begin{align*}
\phi(F,G) & = \psi_0\Big(\bar\xi,F(\bar\xi),\bbH\big(G(\bar\xi)\big)\Big) + h(F,G) + \sum_{k=1}^q \iota_{(-\infty,0]} \bigg(\psi_k\Big(\bar\xi,F(\bar\xi),\bbH\big(G(\bar\xi)\big)\Big)\bigg)\\
& + \iota_C\Big(F(\bar\xi), \bbH\big(G(\bar\xi)\big)\Big) + \iota_{D_\epsilon(\delta)}\big( G(\bar\xi) \big).
\end{align*}
The minimum value of $\phi$ cannot be $\infty$ by assumption (e) and the fact that $h$ is real-valued. Thus,  $(F^\star,G^\star)$ is a minimizer of $\phi$. Moreover, $(F^\star(\bar\xi), \bbH(G^\star(\bar\xi)))\in C$, $G^\star(\bar\xi)\in D_\epsilon(\delta)$, and
\[
\psi_k\Big(\bar\xi,F^\star(\bar\xi),\bbH\big(G^\star(\bar\xi)\big)\Big) \leq 0, ~k=1, \dots, q,
\]
which implies that $(F^\star(\bar\xi),\bbH(G^\star(\bar\xi)))$ is feasible in (AP)$(\bar\xi)$. Let $(x,y)\in C$ satisfy $\psi_k(\bar\xi,x,y)\leq 0$, $k=1, \dots, q$, which then is an arbitrary feasible point in (AP)$(\bar\xi)$. We construct $F$ and $G$ by setting
\begin{align*}
F(\xi) & = F^\star(\xi) - F^\star(\bar\xi) + x\\
g_i(\xi) & = g^\star_i(\xi) - g_i^\star(\bar\xi) + \epsilon(2y_i - 1), ~i=1, \dots, m,
\end{align*}
where $y_i$ is the $i$th component of $y$ and $g_i$ is the $i$th component function of $G$. Note that $\bbH(G(\bar\xi)) = y$. By assumption (c), $F\in \cF$ and $G\in \cG$. Consequently, $(F(\bar\xi), \bbH(G(\bar\xi))) \in C$, $G(\bar\xi) \in D_\epsilon(\delta)$, and
\[
\psi_k\Big(\bar\xi,F(\bar\xi),\bbH\big(G(\bar\xi)\big)\Big)\leq 0, ~k=1, \dots, q,
\]
by construction. Since $(F^\star,G^\star)$ is a minimizer of $\phi$, this implies that
\begin{align*}
& \psi_0\Big(\bar\xi,F^\star(\bar\xi),\bbH\big(G^\star(\bar\xi)\big)\Big) + h(F^\star,G^\star) = \phi(F^\star,G^\star)\\
& \leq \phi(F,G) = \psi_0\Big(\bar\xi,F(\bar\xi),\bbH\big(G(\bar\xi)\big)\Big) + h(F,G)= \psi_0(\bar\xi,x,y) + h(F,G).
\end{align*}
Since $h(F^\star,G^\star) = h(F,G)$ by assumption (d), this implies that
\[
\psi_0\Big(\bar\xi,F^\star(\bar\xi),\bbH\big(G^\star(\bar\xi)\big)\Big) \leq \psi_0(\bar\xi,x,y).
\]
We conclude that $(F^\star(\bar\xi),\bbH(G^\star(\bar\xi)))$ is optimal in (AP)$(\bar\xi)$ because $(x,y)$ is an arbitrary feasible point for that problem. We have already confirmed that $G^\star(\bar \xi)\in D_\epsilon(\delta)$. The fact that $\bbH(\hat G^\nu(\bfxi^\nu(\omega))) = \bbH(G^\star(\bar \xi))$ for all $\omega\in \Omega$ when $\nu\in N$ is sufficiently large follows immediately from $\hat G^\nu(\bfxi^\nu(\omega)) \to G^\star(\bar \xi) \in D_\epsilon(\delta)$ regardless of $\omega\in \Omega$. Since $\phi^\nu(\hat F^\nu,\hat G^\nu)$ converges  along $\nu\in N$ to a value that is not $\infty$, we also conclude that
\begin{align*}
\lim_{\nu \in N} \phi^\nu(\hat F^\nu,\hat G^\nu) & = \lim_{\nu\in N} \cR_0\Big(\psi_0^\nu\big(\bfxi^\nu,\hat F^\nu(\bfxi^\nu),\bbH\big(\hat G^\nu(\bfxi^\nu)\big)\big)\Big) + h^\nu(\hat F^\nu,\hat G^\nu)\\
& = \phi(F^\star,G^\star) = \psi_0\Big(\bar\xi,F^\star(\bar\xi),\bbH\big(G^\star(\bar\xi)\big)\Big) + h(F^\star,G^\star).
\end{align*}

Under the additional assumption that both $\cF$ and $\cG$ contain the zero mappings and $h(F,G)=0$ only when $F$ and $G$ are constant mappings, we realize that $F^\star$ and $G^\star$ must be constant mappings because otherwise one can construct $F^{\star\star}$ and $G^{\star\star}$ by setting
\[
F^{\star\star}(\xi) = F^{\star}(\bar\xi) ~\mbox{ and }~ G^{\star\star}(\xi) = G^{\star}(\bar\xi),
\]
and these mappings are in $\cF$ and $\cG$, respectively. This achieves $\phi(F^{\star\star},G^{\star\star}) < \phi(F^{\star},G^{\star})$ because $0=h(F^{\star\star},G^{\star\star}) < h(F^{\star},G^{\star})$.\eop\\

\state Proof of Corollary \ref{cor:convergence}. The arguments leading to the verification of \eqref{eqn:liminfcond} in the proof of Theorem \ref{thm:convergence} remain valid under the present assumptions. Consequently, using the notation from that proof, we obtain
\[
\infty>\nliminf_{\nu \in N} \gamma^\nu \geq \nliminf_{\nu \in N} \phi^\nu(\hat F^\nu,\hat G^\nu) \geq \phi(F^\star,G^\star) = \psi_0\Big(\bar \xi,F^\star(\bar\xi),\bbH\big(G^\star(\bar\xi)\big)\Big) + h(F^\star,G^\star),
\]
where we use the fact that $\phi(F^\star,G^\star)<\infty$. Thus, $(F^\star(\bar\xi), \bbH(G^\star(\bar\xi)))$ is feasible for (AP)$(\bar\xi)$ and $G^\star(\bar \xi) \in D_\epsilon(\delta)$. The fact that $\bbH(\hat G^\nu(\bfxi^\nu(\omega)))$ $=$ $\bbH(G^\star(\bar \xi))$ for all $\omega\in \Omega$ when $\nu\in N$ is sufficiently large follows because $\hat G^\nu(\bfxi^\nu(\omega)) \to G^\star(\bar \xi) \in D_\epsilon(\delta)$ regardless of $\omega\in \Omega$.\eop\\

\subsection{Implementation Details}\label{subsec:implementationDetails}

In view of the previous study \cite{LejeuneRoysetMa.23}, we solve (SP2)($\xi$) and the various training problems using a linearization approach that leverages the fact that minimizing $\exp(-\alpha z)$ over $z\in \{0, 1, 2, \dots, T\} \cap \cZ$, where $\cZ$ represents constraints, is equivalent to the problem
\begin{equation*}
\nnmin_{z \in \cZ,w_0, \dots, w_T}  \sum_{j=0}^T w_j e^{-j\alpha} \text{ subject to } \sum_{j=1}^T j \, w_j = z, ~\sum_{j=0}^T w_j = 1,~ w_j \in \{0,1\}, ~j =0,1,2, \ldots,T.
\end{equation*}
Implementing this idea, we obtain the following reformulation of (SP2)($\xi$): 
\begin{align*}
\nnmin_{w, y} ~\sum_{i = 1}^I  q_i(\xi) \sum_{j=0}^T w_{i,j,1} \exp\big(-j\alpha(\xi)\big) &\\
  \mbox{subject to }  \sum_{i = 1}^I q_i(\xi) \sum_{j=0}^T w_{i,j,2} \exp\big(-j\alpha(\xi)\big) &\leq \tau\\
  \sum_{c=1}^C y_{c,t} = 1~~~~&\forall t = 1, \dots, T\\
                     \sum_{c'\in N(c)} y_{c',t-1} \geq y_{c,t} ~~~~&\forall c = 1, \dots, C, ~t=2, \dots, T\\            
                     \sum_{j=1}^T j \, w_{i,j,k} = \sum_{c = 1}^C \sum_{t = 1}^T \zeta_{c,t,i}^k \, y_{c,t}~~~~&\forall i=1, \dots, I, ~k = 1, 2\\
                     \sum_{j=0}^T w_{i,j,k} = 1 ~~~~&\forall i = 1, \dots, I, ~k = 1, 2\\
                      w_{i,j,k}, y_{c,t}  \in \{0,1\}~~~~\forall c=1, \dots, C,~ t=1, \dots, T, ~&i=1, \dots, I, ~j=0, \dots, T, ~k=1, 2.
\end{align*}
The variables $w_{i,j,k}$ can be relaxed to continuous on $[0,1]$ without jeopardizing the equivalence with (SP2)($\xi$) because the exponential function is strictly convex. This relaxation, however, tends to cause numerical instabilities and we use the stated formulation consistently. We linearize (SP1)($\xi$) and the training problems (EW-SP2)$^\nu$, (WW-SP2)$^\nu$, and ($\beta$-SP1)$^\nu$ similarly.\\

Under the expectation risk measure for $\cR_0$, the training problem (TP)$^\nu$ supporting (SP2)($\xi$) takes the following specialized form after leveraging a similar expansion of variables as in (TP-super): 
\begin{align}
\mbox{(EW-SP2)}^\nu~~~\nnmin_{y,B,b} \frac{1}{|\Omega|}\sum_{\omega \in \Omega} \sum_{i = 1}^I  q_i\big(\bfxi^\nu(\omega)\big) \exp\bigg(& - \alpha\big(\bfxi^\nu(\omega)\big) \sum_{c = 1}^C \sum_{t = 1}^T \zeta_{c,t,i}^1 \, y_{c,t}(\omega)  \bigg) + \theta \sum_{c=1}^C \sum_{t=1}^T \|B_{c,t}\|_1\nonumber\\
  \mbox{subject to } \sum_{i = 1}^I q_i\big(\bfxi^\nu(\omega)\big) \exp\bigg(- \alpha\big(\bfxi^\nu(\omega)\big) \sum_{c = 1}^C &\sum_{t = 1}^T \zeta_{c,t,i}^2 \, y_{c,t}(\omega)  \bigg) \leq \tau ~~\forall \omega \in \Omega\label{eqn:vfirstconSPnu}\\
  \sum_{c=1}^C y_{c,t}(\omega)  = 1~~~~&\forall t = 1, \dots, T, ~\omega \in \Omega\label{eqn:vSecconSPnu}\\
                     \sum_{c'\in N(c)} y_{c',t-1}(\omega)  \geq y_{c,t}(\omega) ~~~~&\forall c = 1, \dots, C, ~t=2, \dots, T, ~\omega \in \Omega\label{eqn:vthridconSPnu}\\
                     -\delta + (\delta + \epsilon) y_{c,t}(\omega)  \leq \big\langle B_{c,t}, \bfxi^\nu(\omega)\big\rangle + b_{c,t} ~~~~ &\forall c=1, \dots, C, ~t=1, \dots, T, ~\omega \in \Omega\label{eqn:firstconSPnu}\\
                     - \epsilon + (\delta + \epsilon)y_{c,t}(\omega)  \geq \big\langle B_{c,t}, \bfxi^\nu(\omega)\big\rangle + b_{c,t} ~~~~ &\forall c=1, \dots, C, ~t=1, \dots, T, ~\omega \in \Omega\label{eqn:secconSPnu}\\
                     B_{c,t}  \in \reals^{r}, ~b_{c,t}  \in \reals~~~~&\forall c=1, \dots, C, ~t=1, \dots, T\label{eqn:vfourconSPnu}\\
                     y_{c,t}(\omega) \in \{0,1\}~~~~&\forall c=1, \dots, C, ~t=1, \dots, T, ~\omega \in \Omega.\label{eqn:vsixconSPnu}
\end{align}

Under the worst-case risk measure for $\cR_0$, the training problem takes the form:
\begin{align*}
\mbox{(WW-SP2)}^\nu~~~~~~~\nnmin_{y,B,b,\gamma} \,\gamma + \theta \sum_{c=1}^C \sum_{t=1}^T \|B_{c,t}\|_1&\nonumber\\
  \mbox{subject to }  ~~\eqref{eqn:vfirstconSPnu}\mbox{-}\eqref{eqn:vsixconSPnu},~~~~~~~~~~~~\gamma & \in \reals\nonumber\\
                        \sum_{i = 1}^I q_i\big(\bfxi^\nu(\omega)\big) \exp\bigg( - \alpha\big(\bfxi^\nu(\omega)\big) \sum_{c = 1}^C \sum_{t = 1}^T \zeta_{c,t,i}^1 \, y_{c,t}(\omega)  \bigg) &\leq \gamma ~~\forall \omega \in \Omega.\nonumber                        
\end{align*}

Under the $\beta$-superquantile risk measure for $\cR_0$, the training problem (TP)$^\nu$ supporting (SP1)($\xi$) takes the form:
\begin{align*}
\mbox{($\beta$-SP1)}^\nu~~~~~~~\nnmin_{y,B,b,\gamma,u} \,\gamma + \frac{1}{1-\beta}\frac{1}{|\Omega|}\sum_{\omega\in \Omega}u(\omega) + \theta \sum_{c=1}^C \sum_{t=1}^T \|B_{c,t}\|_1&\nonumber\\
  \mbox{subject to }  ~~\eqref{eqn:vSecconSPnu}\mbox{-}\eqref{eqn:vsixconSPnu},~~~\gamma \in \reals,  ~~~u(\omega)&\geq 0 ~~\forall \omega \in \Omega.\nonumber\\
                        \sum_{i = 1}^I q_i\big(\bfxi^\nu(\omega)\big) \exp\bigg( - \alpha\big(\bfxi^\nu(\omega)\big) \sum_{c = 1}^C \sum_{t = 1}^T \zeta_{c,t,i}^1 \, y_{c,t}(\omega)  \bigg) - \gamma &\leq u(\omega) ~~\forall \omega \in \Omega.\nonumber                       
\end{align*}

\subsection{Additional Numerical Results}\label{subsub:addnum}

Figure \ref{fig:grid} illustrates the discretized environment consisting of a 9-by-9 grid, with cell numbers, as utilized in Section \ref{sec:num}.

\drawing{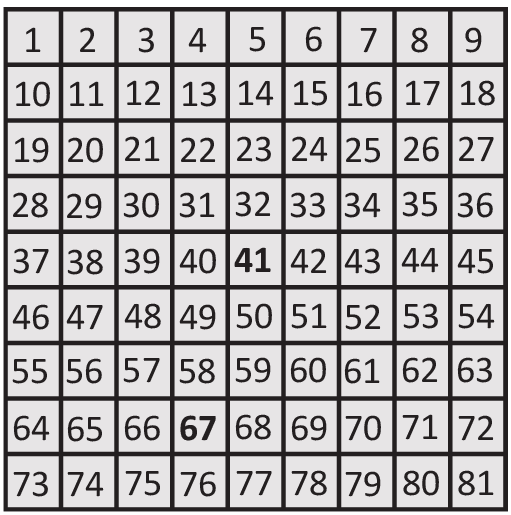}{1.7in} {Cell numbers in discretized environment with 81 cells.}{fig:grid}

Table \ref{runtimes} summarizes computing times for the Decomposition Algorithm across 16 instances of (EW-SP2)$^\nu$ and 16 instances of (WW-SP2)$^\nu$.

\begin{table}[ht]
		\centering
		\setlength\extrarowheight{1pt}
		\footnotesize{
			\begin{tabular}{c|c|c|c|c}
				\hline
				\multirow{2}{*}{$\nu$} 
				&\multicolumn{2}{c|}{(EW-SP2)$^\nu$  } 				&\multicolumn{2}{c}{(WW-SP2)$^\nu$}\\ 
				\cline{2-5} 
				& $|\Omega| = 10$ & $|\Omega| = 100$ & $|\Omega| = 10$ &  $|\Omega| = 100$\\
				\hline
				$1$      & 59  &831 & 39 & 902\\ 			
				$2$      & 330 &786 & 26 & 860\\		
				$3$      & 99  &780 & 32 & 902\\ 
				$4$      & 26  &843 & 19 & 834\\
				$5$      & 32  &841 & 22 & 1099\\ 
				$6$      & 55  &879 & 34 & 1133\\ 
				$7$      & 39  &932 & 37 & 1128\\ 
				$8$      & 38  &952 & 33 & 978\\ 
				\hline
			\end{tabular}
			\caption{Run times (sec.) for the Decomposition Algorithm on different training problems.}\label{runtimes}
		}
\end{table}

Tables \ref{tab:testSextra} and \ref{tab:testSdifferentextra} supplement Tables \ref{tab:testS} and \ref{tab:testSdifferent} by considering the uniform test data for decision rules based on the beta training data.

\begin{table}[ht]
	\centering
	\setlength\extrarowheight{1pt}
	\footnotesize{
		\begin{tabular}{c|c|c|c|c|c|c}
			\hline
			\multirow{2}{*}{rule} 
			& training & test	&number 
			&\multicolumn{3}{c}{suboptimality}\\ 
			\cline{5-7} 
			& data & data & feasible          & min  &avg   & max\\
			\hline
			$B^\nu,b^\nu$ & beta & unif & 38  &0.000 &0.005 &0.017\\	
			MDR           & beta & unif & 100 &0.006 &0.017 &0.031\\	
			AMDR          & beta & unif & 100 &0.000 &0.000 &0.002\\				
			\hline
		\end{tabular}
		\caption{Performance of decision rules obtain from ($\beta$-SP1)$^\nu$. Suboptimality statistics are computed over feasible decisions.}
		\label{tab:testSextra}
	}	
\end{table}

\begin{table}[ht]
	\centering
	\setlength\extrarowheight{1pt}
	\footnotesize{
		\begin{tabular}{c|c|c|c|c|c|c}
			\hline
			\multirow{2}{*}{rule} 
			& training & test	&number 
			&\multicolumn{3}{c}{suboptimality}\\ 
			\cline{5-7} 
			& data & data & feasible          & min  &avg   & max\\
			\hline
			$B^\nu,b^\nu$ & beta & unif & 33  &0.010 &0.048 &0.076\\	
			MDR           & beta & unif & 100 &0.033 &0.052 &0.095\\	
			AMDR          & beta & unif & 100 &0.000 &0.000 &0.007\\				
        \hline
		\end{tabular}
		\caption{Performance of decision rules obtain from ($\beta$-SP1)$^\nu$ against dispersed target. Suboptimality statistics are computed over feasible decisions.}
		\label{tab:testSdifferentextra}
	}	
\end{table}

We supplement the numerical results of Subsection \ref{subsec:feasRules}, which utilizes superquantiles for risk measures in the training problem ($\beta$-SP1)$^\nu$, by also considering expectations and worst-case risk. Let (E-SP1)$^\nu$ and (W-SP1)$^\nu$ be the training problems obtained from (EW-SP2)$^\nu$ and (WW-SP2)$^\nu$, respectively, by dropping \eqref{eqn:vfirstconSPnu} pertaining to the second target. 

Using the uniform training data, the Decomposition Algorithm solves (E-SP1)$^\nu$ in 335 seconds to a relative optimality gap of 0.06\%, which corresponds to an absolute gap of 0.002. This produces a decision rule $(B^\nu,b^\nu)$. The first six rows of Table \ref{tab:testE} report the performance of that decision rule as well as its derivatives MDR and AMDR. 
We find that the three decision rules perform equally well. This is caused by the fact that $B^\nu$ actually is the zero matrix in this case. Thus, for the range of parameter values considered in training, no adaptation is found to be beneficial. 

For the beta training data, the Decomposition Algorithm solves (E-SP1)$^\nu$ in 936 seconds, with relative optimality gap of 4.12\% and absolute gap of 0.015. The last six rows of Table \ref{tab:testE} report the performance of the resulting decision rules. Now, direct use of  $(B^\nu,b^\nu)$ is feasible only at a third of the test points. Moreover, the importance of adaptation is clear: AMDR outperforms MDR with a wide margin.


\begin{table}[ht]
	\centering
	\setlength\extrarowheight{1pt}
	\footnotesize{
		\begin{tabular}{c|c|c|c|c|c|c}
			\hline
			\multirow{2}{*}{rule} 
			& training & test	&number 
			&\multicolumn{3}{c}{suboptimality}\\ 
			\cline{5-7} 
			& data & data & feasible          & min  &avg   & max\\
			\hline
			$B^\nu,b^\nu$ & unif & unif & 100 &0.000 &0.002 &	0.013\\	
			MDR           & unif & unif & 100 &0.000 &0.002 &	0.013\\	
			AMDR          & unif & unif & 100 &0.000 &0.002 &	0.013\\				
			\hline
			$B^\nu,b^\nu$ & unif & beta & 100 &0.000 &0.009 &	0.033\\	
			MDR           & unif & beta & 100 &0.000 &0.009 &	0.033\\	
			AMDR          & unif & beta & 100 &0.000 &0.009 &	0.033\\				
			\hline
			$B^\nu,b^\nu$ & beta & unif & 30  &0.000 &0.002 &	0.013\\	
			MDR           & beta & unif & 100 &0.000 &0.011 &	0.026\\	
			AMDR          & beta & unif & 100 &0.000 &0.000 &	0.006\\				
			\hline
			$B^\nu,b^\nu$ & beta & beta & 37  &0.000 &0.009 &	0.033\\	
			MDR           & beta & beta & 100 &0.001 &0.017 &	0.049\\	
			AMDR          & beta & beta & 100 &0.000 &0.002 &	0.014\\
        \hline
		\end{tabular}
		\caption{Performance of decision rules obtain from (E-SP1)$^\nu$. Suboptimality statistics are computed over feasible decisions.}
		\label{tab:testE}
	}	
\end{table}

Table \ref{tab:testW} provides parallel results to those in Table \ref{tab:testE} after the switch from (E-SP1)$^\nu$ to (W-SP1)$^\nu$. The solution of (W-SP1)$^\nu$ using the Decomposition Algorithm takes 361 seconds on the uniform training data, producing a relative optimality gap of 2.25\% (absolute gap of 0.009). On the beta training data, the numbers are 511 seconds, 1.89\%, and 0.008. (W-SP1)$^\nu$ produces mostly worst levels of suboptimality compare to those obtained from (E-SP1)$^\nu$ for the instances in Table \ref{tab:testW}, but does achieve more feasible decisions.  Still, AMDR produced by (W-SP1)$^\nu$ consistently prescribes quality decisions.


\begin{table}[ht]
	\centering
	\setlength\extrarowheight{1pt}
	\footnotesize{
		\begin{tabular}{c|c|c|c|c|c|c}
			\hline
			\multirow{2}{*}{rule} 
			& training & test	&number 
			&\multicolumn{3}{c}{suboptimality}\\ 
			\cline{5-7} 
			& data & data & feasible          & min  &avg   & max\\
			\hline
			$B^\nu,b^\nu$ & unif & unif &  92 &0.000 &0.009 &0.023\\	
			MDR           & unif & unif & 100 &0.003 &0.011 &0.023\\	
			AMDR          & unif & unif & 100 &0.000 &0.002 &0.013\\				
			\hline
			$B^\nu,b^\nu$ & unif & beta & 70  &0.000 &0.016 &0.034\\	
			MDR           & unif & beta & 100 &0.001 &0.018 &0.042\\	
			AMDR          & unif & beta & 100 &0.000 &0.008 &0.033\\				
			\hline
			$B^\nu,b^\nu$ & beta & unif & 67  &0.000 &0.008 &0.018\\	
			MDR           & beta & unif & 100 &0.003 &0.011 &0.023\\	
			AMDR          & beta & unif & 100 &0.000 &0.001 &0.013\\				
			\hline
			$B^\nu,b^\nu$ & beta & beta & 65  &0.000 &0.015 &0.032\\	
			MDR           & beta & beta & 100 &0.001 &0.018 &0.042\\	
			AMDR          & beta & beta & 100 &0.000 &0.006 &0.023\\
        \hline
		\end{tabular}
		\caption{Performance of decision rules obtain from (W-SP1)$^\nu$. Suboptimality statistics are computed over feasible decisions.}
		\label{tab:testW}
	}	
\end{table}

We replicate the results of Tables \ref{tab:testE} and \ref{tab:testW} also for the more challenging dispersed target case. Now, the solution of (E-SP1)$^\nu$ using the Decomposition Algorithm takes 746 seconds on the uniform training data, producing a relative optimality gap of 2.05\% (absolute gap of 0.011). On the beta training data, the numbers are 130 seconds, 4.70\%, and 0.026. For (W-SP1)$^\nu$, the Decomposition Algorithm takes 940 seconds on the uniform training data, producing a relative optimality gap of 2.17\% (absolute gap of 0.022). On the beta training data, the numbers are 830 seconds, 3.97\%, and 0.024. Tables \ref{tab:testEdifferent} and \ref{tab:testWdifferent} report on the performance of the resulting decision rules. While the level of suboptimality of the decision rules are higher compared to the original target setting, AMDR remains viable. Again, (E-SP1)$^\nu$ tends to produce better levels of suboptimality (see Table \ref{tab:testEdifferent}) as compared to (W-SP1)$^\nu$ in Table \ref{tab:testWdifferent}, but is also more likely to produce infeasible decisions. 


\begin{table}[ht]
	\centering
	\setlength\extrarowheight{1pt}
	\footnotesize{
		\begin{tabular}{c|c|c|c|c|c|c}
			\hline
			\multirow{2}{*}{rule} 
			& training & test	&number 
			&\multicolumn{3}{c}{suboptimality}\\ 
			\cline{5-7} 
			& data & data & feasible            & min  &avg   & max\\
			\hline
			$B^\nu,b^\nu$ & unif & unif & 100   &0.000 &0.011 &0.032\\	
			MDR           & unif & unif & 100   &0.000 &0.011 &0.032\\	
			AMDR          & unif & unif & 100   &0.000 &0.011 &0.032\\				
			\hline
			$B^\nu,b^\nu$ & unif & beta & 100   &0.000 &0.028 &0.086\\	
			MDR           & unif & beta & 100   &0.000 &0.028 &0.086\\	
			AMDR          & unif & beta & 100   &0.000 &0.028 &0.086\\				
			\hline
			$B^\nu,b^\nu$ & beta & unif & 1     &0.028 &0.028 &0.028\\	
			MDR           & beta & unif & 100   &0.000 &0.018 &0.075\\	
			AMDR          & beta & unif & 100   &0.000 &0.001 &0.021\\				
			\hline
			$B^\nu,b^\nu$ & beta & beta &  5    &0.000 &0.023 &0.094\\	
			MDR           & beta & beta & 100   &0.000 &0.018 &0.075\\	
			AMDR          & beta & beta & 100   &0.000 &0.001 &0.021\\
        \hline
		\end{tabular}
		\caption{Performance of decision rules obtain from (E-SP1)$^\nu$ against dispersed target. Suboptimality statistics are computed over feasible decisions.}
		\label{tab:testEdifferent}
	}	
\end{table}


\begin{table}[ht]
	\centering
	\setlength\extrarowheight{1pt}
	\footnotesize{
		\begin{tabular}{c|c|c|c|c|c|c}
			\hline
			\multirow{2}{*}{rule} 
			& training & test	&number 
			&\multicolumn{3}{c}{suboptimality}\\ 
			\cline{5-7} 
			& data & data & feasible          & min  &avg   & max\\
			\hline
			$B^\nu,b^\nu$ & unif & unif & 100 &0.006 &0.023 &0.049\\	
			MDR           & unif & unif & 100 &0.009 &0.024 &0.049\\	
			AMDR          & unif & unif & 100 &0.000 &0.003 &0.027\\				
			\hline
			$B^\nu,b^\nu$ & unif & beta & 100 &0.000 &0.036 &0.122\\	
			MDR           & unif & beta & 100 &0.005 &0.038 &0.094\\	
			AMDR          & unif & beta & 100 &0.000 &0.016 &0.068\\				
			\hline
			$B^\nu,b^\nu$ & beta & unif & 50  &0.000 &0.006 &0.044\\	
			MDR           & beta & unif & 100 &0.008 &0.018 &0.036\\	
			AMDR          & beta & unif & 100 &0.000 &0.002 &0.025\\				
			\hline
			$B^\nu,b^\nu$ & beta & beta &  36 &0.000 &0.027 &0.084\\	
			MDR           & beta & beta & 100 &0.002 &0.032 &0.080\\	
			AMDR          & beta & beta & 100 &0.000 &0.013 &0.068\\
        \hline
		\end{tabular}
		\caption{Performance of decision rules obtain from (W-SP1)$^\nu$ against dispersed target. Suboptimality statistics are computed over feasible decisions.}
		\label{tab:testWdifferent}
	}	
\end{table}

\end{document}